\documentclass{book}

\usepackage{amsmath, amssymb, amsfonts, amsthm}
\usepackage{graphicx}
\usepackage{natbib}
\usepackage{color}
\usepackage{bm}
\usepackage{subfigure}
\usepackage{graphicx}
\usepackage{mathabx}
\usepackage{multirow}
\usepackage{setspace}

\usepackage{ResCom}
\usepackage{natbib}

\author{Hyeongkwan Kim}
\title{Gersten-Witt Complex of Hirzebruch Surfaces}

\date{2013} 

\usepackage[pdftex, pdfusetitle, plainpages=false, 
				letterpaper, bookmarks, bookmarksnumbered,
				colorlinks, linkcolor=black, citecolor=black,
	         filecolor=black, urlcolor=black]
				{hyperref}

\begin{document}
\maketitle

\chapter{Introduction}

We will first review some basic definitions of  symmetric bilinear forms, and define Witt groups of fields, rings, and schemes. In Section \ref{GWC}, we introduce the notion of the Gersten-Witt complex of a scheme. We will see that for certain class of schemes with nice geometric properties, the Gersten-Witt complex is the global section of a flasque resolution of a Witt sheaf on the scheme. We will see how its cohomologies can be computed easily if the scheme is a complex toric variety.

\section{Symmetric bilinear forms over a field}

Let $k$ be a field with $\chr k\neq2$. A \itk{symmetric bilinear form} over $k$ is a map of the form
$$\phi:V\times V\ra k,$$
where $V$ is a finite-dimensional $k$-\vsp, \st
\bea{\phi(v,w) &=& \phi(w,v), \\
\phi(v+v',w) &=& \phi(v,w)+\phi(v',w), \\
\phi(a v,w) &=& a\phi(v,w),
} for every $a\in k$ and $v,v',w\in V$.
We will denote the form by $(V,\phi)$. It is \itk{\ns} if the induced map
$$\ad\phi:V\ra\Hom(V,k)$$
is bijective, and \itk{anisotropic} if $\phi(v,v)=0$ implies $v=0$.

Two bilinear forms $(V,\phi)$ and $(W,\psi)$  are \itk{isometric} if there is an \iso\ of vector spaces $f:V\ras W$ \st the diagram
\xym{V\times V \ar[r]^\phi \ar[d]_(.45){f\times f}^(.45){\rotatebox{90}{$\sim$}} & k \\
W\times W \ar[ru]_(.55)\psi
} commutes.

If $W\sst V$ is a subspace, we define its \itk{orthogonal complement}
$$W\pp:=\setst{v\in V}{\phi(v,w)\fall w\in W}.$$
If $W\sst W\pp$, then $W$ is \itk{totally isotropic}, and if $W=W\pp$, then $W$ is \itk{orthogonal}. There is a dimension equation \cite[1.3]{lam},
\beq{\dm W+\dm W\pp=\dm V.\label{dim}}
A \ns\ form $(V,\phi)$ is \itk{hyperbolic} if $V$ has an orthogonal subspace, or equivalently, a totally isotropic subspace of half the dimension (by (\ref{dim})).

Every symmetric bilinear form over $k$ can be represented by a symmetric matrix $M$, and it is \ns\ \ifof $\det M\neq0$. The forms represented by matrices $M$ and $M'$ are isometric \ifof there is a \ns\ matrix $Q$ \st $M'=QMQ^T$. Every symmetric bilinear form over $k$ can be diagonalized by such a transformation.

The hyperbolic form of dimension 2 is called the \itk{hyperbolic plane}. It is represented by a matrix
$$\left(\begin{array}{cc} 0 & 1 \\ 1 & 0\end{array}\right),$$
which diagonalizes to (recall our assumption that $2\in A$ is a unit)
$$\left(\begin{array}{cc} 1 & 0 \\ 0 & -1 \end{array}\right)=\left(\begin{array}{cc} 1 & 1/2 \\ 1 & -1/2 \end{array}\right)\left(\begin{array}{cc} 0 & 1 \\ 1 & 0 \end{array}\right)\left(\begin{array}{cc} 1 & 1 \\ 1/2 & -1/2 \end{array}\right).$$
It can be shown that every hyperbolic space decomposes into a direct sum of hyperbolic planes \cite[3.4(1)]{lam}.

\section{Witt groups}

\subsection{Witt group of a field}

Let $k$ be as defined above, and let $Q(k)$ be the set of isometry classes of \ns\ symmetric bilinear forms over $k$. $Q(k)$ is a semigroup, where the addition is defined by the orthogonal sum,
$$[V,\phi]+[W,\psi]=[V\opl W,\phi\opl\psi].$$
The Grothendieck group of  $Q(k)$ modulo the subgroup generated by hyperbolic forms is called the \itk{Witt group of $k$}, denoted by $W(k)$. By the diagonalizability, every element of $W(k)$ can be represented by a finite sum of unary forms
$$\ang{a_1}+\ang{a_2}+\cdots+\ang{a_r},$$
where $a_1,\ldots,a_r\in k\cro$. Quotienting by hyperbolic forms allows one to write
$$\ang{-a}=-\ang{a}\in W(k).$$
Furthermore, it is a consequence of Witt decomposition theorem \cite[4.1]{lam} that every element of $W(k)$ can be represented by an anisotropic form.

Note that the isometry relation implies  $\ang{a^2}=\ang{1}$.

Let us look at some examples. The signature map and the dimension map respectivly induce \isos\ \cite[p.~41-42]{lam}
$$W(\re)\ras\Z,\qqq W(\cx)\ras\Z/2.$$
If $p\in\Z$ is an odd prime, the Witt group of the finite field $\Ff_p$ is given by\footnote{$W(k)$ admits a ring structure induced by tensor product, but we won't need this in this paper.} \cite[p.~45]{lam}
$$W(\Ff_p)=\begin{cases}
            (\Z/2\Z)^2 & \quad\mb{ if }\quad p\equiv 1\  (\operatorname{mod} 4), \\
\Z/4\Z & \quad\mb{ if }\quad p\equiv 3\  (\operatorname{mod} 4).
           \end{cases}
$$
The Witt group of the rationals is given by \cite[p.~175]{lam}
\beq{W(\Q)\simeq\Z\opl\Z/2\opl\coprod_{p\neq2}W(\Z/p\Z).\label{hasse}}
We will see that this is derived from the Gersten-Witt complex of $\spec\Z$ (\ref{iwill}).

We will come back to this later, but for now it suffices to say that the Gersten-Witt complex is largely an attemp to generalize this \iso\ to schemes.

For future reference, we state a theorem by Springer and Knebusch \cite[p.~85]{milnor}:

\thm{Let $A$ be a discrete valuation ring with maximal ideal $\maxi$ and quotient field $F$, where $\chr (A/\maxi)\neq2$. If $\pi\in A$ is a generator of $\maxi$ and $i=1$ or $2$,  there is a unique \hmm
$$\pl^\pi_i:W(F)\ra W(A/\maxi)$$
\st
$$\ang{\pi^ju}\mt\begin{cases}
                  \ang{\bar{u}} & \quad\mb{ if }\quad j\not\equiv i \pmod{2},\\
0 & \quad\mb{ if }\quad j\equiv i \pmod{2}.
                 \end{cases}
$$
\label{springer}
}
Note that $\pl_2^\pi$ depends on the choice of the generator $\pi$, while $\pl_1^\pi$ doesn't. $\pl_1^\pi$ and $\pl_2^\pi$ are called the \itk{first and second residue \hmm}, \resp.

\subsection{Witt group of a ring}

The Witt group can be similarly defined for a ring $A$ in which 2 is a unit.  The finite-dimensional $k$-vector spaces are replaced by finitely generated projective \amods, the rank of projective modules replacing the dimension of vector spaces. The definition of nonsingularity remains the same (i.e.,  bijectivity of the adjoint).
Instead of quotienting the Grothendieck group associated with the semigroup of \ns\ symmetric bilinear forms by hyperbolic forms, we quotient by a larger class of forms called lagrangians\footnote{Knebusch\cite{Kn2} uses the term ``metabolic space" for our lagrangian, ``split metabolic space" for our hyperbolic form, ``lagrangian" for our sublagrangian, and ``sublagrangian" for our totally isotropic space.}: if $M$ is a finitely generated projective \amod, a \ns\ symmetric bilinear form
$$\phi:M\times M\ra A$$
is called a \itk{lagrangian} if there is a direct summand $N\sst M$ \st $\phi|_{N\times N}=0$ and the induced pairing
$$N\times(M/N)\ra A$$
is \ns\ (i.e., both adjoints are bijective\footnote{In fact, the reflexivity of finitely generated projective modules implies that one adjoint is bijective \ifof the other is.}). The submodule $N$ is called a \itk{sublagrangian}. As in the case of hyperbolic spaces, $(M,\phi)$ is a lagrangian \ifof $M$ has  an orthogonal direct summand, or equivalently, a totally isotropic direct summand with half the rank of $M$ \cite[Corollary 2, ii]{Kn2}.

The hyperbolic form still plays a role, mainly due to its useful properties (e.g., (\ref{useful}) below). Let us first generalize hyperbolic forms in the context of finitely generated projective modules. If $M$ is a finitely generated projective \amod\ and $M\du:=\Hom(M,A)$,  the \itk{hyperbolic form}  associated with $M$ is defined to be a symmetric bilinear form
$$\phi:(M\opl M)\du\times (M\opl M\du)\ra A$$
induced by the canonical pairing $M\times M\du\ra A$, and requiring that 
$$\phi|_{M\times M}=0,\qqq \phi|_{ M\du\times M\du}=0.$$
The reflexivity of finitely generated projective modules ensures its nonsingularity. Note that if $\rho_M:M\ra M\ddu$ is the canonical map, $(M,\phi)$ can be represented by the ``matrix"
\beq{\left(\begin{array}{cc} 0 & \id_{M\du} \\ \rho_M & 0\end{array}\right).\label{hyper}}
$\rho_M$ is an \iso\ because $M$ is finitely generated projective. Hence, if $M$ is free, then (\ref{hyper}) is isometric to a direct sum of the hyperbolic planes,
$$\left(\begin{array}{cc} 0 & 1 \\ 1 & 0\end{array}\right).$$
Thus our definition of hyperbolic form agrees with the previous one over fields. Note that the submodule $M\sst M\opl M\du$ is a sublagrangian, since the canonical pairing
$$M\times M\du\ra V$$
is \ns\ by the reflexivity of finitely generated projective modules. Hence, every hyperbolic form is a lagrangian.

We note a useful lemma that we will be needed later:

\lemm{If $(M,\phi)$ is a \ns\ symmetric bilinear form, then the form $(M\opl M,\phi\opl(-\phi))$ is isometric to the hyperbolic form associated with $M$.
\label{useful}
}

\pf{Let
$$\Delta:=\setst{(m,m)\in M\opl M}{m\in M},\qq N:=\setst{(m/2,-m/2)\in M\opl M}{m\in M}.$$
Clearly $\Delta, N\simeq M$, $\Delta\cap N=\set{0}$, and the \iso
$$M\opl M\ras\Delta+N,\qqq(m,m')\mt(m+m/2,m'-m'/2)$$
establishes an isometry between $(M\opl M,\phi\opl(-\phi))\simeq(\Delta\opl N,\phi\opl(-\phi))$ and the hyperbolic form associated with $M$.
}

Let us look at some examples.
The signature map induces an \iso\ $W(\Z)\ras\Z$ \cite[p.~23]{milnor}, and  there is a split short exact sequence \cite[p.~175]{lam}
\beq{0\ra W(\Z)\ra W(\Q)\ra (\Z/2)\opl\coprod_{p\neq2} W(\Z/p\Z)\ra0,\label{iwill}}
which gives rise to the \iso\ (\ref{hasse}). The split exactness of (\ref{iwill}) is a consequence of the Hasse-Minkowski principle applied to the global field $\Q$. More generally, if we use a Gorenstein ring $A$ of dimension $n$ instead of $\Z$, we can construct a complex of the form
\beq{0\ra W(A)\ra \coprod_{\hei\pr=0} W(\kappa(\pr))\ra\cdots\ra\coprod_{\hei\pr=n} W(\kappa(\pr))\ra0,\label{afterd}}
where $\kappa(\pr)$ is the residue class field at $\pr\in\spec A$. However, (\ref{afterd}) is not exact in general. As we shall see, (\ref{afterd}) is the Gersten-Witt complex of $\spec A$, and Pardon  \cite[5.1]{bigP} proved that it is exact if $A$ is a regular local ring and is of essentially finite type over a field of characteristic different from 2.

\subsection{Witt group of a scheme}

Knebusch\cite{Kn2} defined the Witt group of a scheme $(X,\oox)$, whose elements are represented by symmetric bilinear forms
$$\phi:\M\times\M\ra\oo_X,$$
where $\M$ is a locally free sheaf of $\oox$-modules. The definition of nonsingularity remains the same as in the affine case (i.e., both adjoints are bijective), but the sublagrangian is no longer required to be a split submodule. In fact, sublagrangians always split for affine schemes \cite[p.~134]{Kn2}, so this definition agrees with the earlier one. The criteria for a sublagrangian is thus an orthogonal submodule, or equivalently, a totally isotropic submodule with half the rank  \cite[Corollary 2ii]{Kn2}. Moreover, the rank of a totally isotropic submodule cannot exceed half the rank \cite[Corollary 2i]{Kn2}.

\section{Gersten-Witt complex}
\label{GWC}

The \itk{Gersten-Witt complex} of a scheme $X$ of dimension $n$ is a complex of the form
$$0\ra \coprod_{x\in X^{(0)}}W(\kappa(x))\rf{d^0}\coprod_{x\in X^{(1)}}W(\kappa(x))\rf{d^1}\cdots\ra\coprod_{x\in X^{(n)}}W(\kappa(x))\rf{d^n}0,$$
where $\kx$ is the residue class field at $x\in X$, and $X^{(p)}\sst X$ is the subset of codimension $p$ points.
A main challenge in constructing such a complex is to define canonical boundary maps $d^p$. The second residue \hmm\ (Theorem \ref{springer}) depends on the choice of the uniformizer, so it cannot be used directly. Several authors constructed the complex using different methods, such as spectral sequences~\cite{ballmer}, the module of differentials~\cite{knebusch,schmid}, and the canonical sheaf~\cite{bigP}. In the latter, using the analogues of Quillen's localization sequence of $K$-groups \cite[8.4]{swan}, Pardon constructed the Gersten-Witt complex of a Gorenstein scheme, and showed that it is acyclic if the scheme is the  spectrum of a regular local ring essentially of finite type over a field with characteristic different from 2. This implies the existence of a flasque resolution $\W\sbu(X)$ of a \itk{Witt sheaf} of a regular scheme $X$ of finite type over $k$, whose cohomologies  furnish us with a new set of invariants for the scheme $X$.

Computing cohomologies of the Gersten-Witt complex is difficult in practice, because it involves Witt groups of residue class fields at all (possibly infinitely many) points of the scheme. 
On the other hand, if $X$ is a complex $n$-dimensional toric variety, then $X$ is filtered as
$$X=X^0\supset X^1\supset\cdots\supset X^n\supset X^{n+1}=\ets,$$
where $X^p-X^{p+1}$ is a disjoint union of $(n-p)$-tori, $\spec\cx[x_1,1/x_1,\ldots,x_{n-p},1/x_{n-p}]$.
Takeda~\cite{takeda} showed the Gersten-Witt complex of $K$-groups is quasi-\isoc\ to a complex of $K$-groups of coherent sheaves of tori. We will show that the same result holds for Witt groups. The Witt group of the $n$-torus is known; for example  \cite{kar, ran},
\begin{eqnarray}
W(\cx[x,1/x]) &=& (\Z/2)^2,\label{ff1}\\
W(\cx[x,y,1/x,1/y])&=& (\Z/2)^4.\label{ff2}
\end{eqnarray}
Hence, the quasi-\isoc\ complex would consist of a finite number of Witt groups which are finite-dimensional vector spaces, and its cohomologies are much easier to compute. Using this method, we will compute cohomologies of the Gersten-Witt complex of the  toric variety Hirzebruch surface $\hn$. Specifically, the Gersten-Witt complex of $\hn$ is given by
\beq{0\ra \coprod_{x\in H_n^{(0)}}W(\kappa(x))\rf{d^0}\coprod_{x\in H_n^{(1)}}W(\kappa(x))\rf{d^1}\coprod_{x\in H_n^{(2)}}W(\kappa(x))\rf{d^2}0,\label{GW}}
and we will show that it is quasi-\isoc\ to a complex of the form
\beq{0\ra W(H_n-H_n^1)\rf{d^0_{\hn}} W(H_n^1-H_n^2)\rf{d^1_{\hn}} W(H_n^2)\ra0,\label{toricc}}
where $H^p_n$ is the closure of codimension $p$ orbits of the torus action,
$$H_n=H_n^0\supset H_n^1\supset H_n^2\supset H_n^3=\ets,$$
and $\hn^p-\hn^{p+1}$ is a finite union of $(n-p)$-tori.

The rest of this paper is organized as follows:

In Chapter \ref{coeff}, we will review Pardon's construction of the Gersten-Witt complex of Gorenstein schemes.

In Chapter \ref{hirze}, we will introduce the Hirzebruch surface $\hn$, and construct a toric complex which is quasi-\isoc\ to the Gersten-Witt complex of $\hn$.

In Chapter \ref{quasi}, we will prove the quasi-\iso.

In Chapter \ref{compu}, we will compute the boundary maps of the toric complex.

In Chapter \ref{coho}, we will compute cohomologies of the toric complex, and find that 
$$H^0(\W\sbu(\hn))= H^0(\W\sbu(\hn))=\Z/2\ffall n\in\Z,$$
but
$$H^i(\W\sbu(H_{\mb{\scriptsize even}}))\neq H^i(\W\sbu(H_{\mb{\scriptsize odd}}))\qq\mb{for}\quad i=1,2.$$

\chapter{Witt groups with coefficients}
\label{coeff}

To construct a canonical Gersten-Witt complex of a Gorenstein scheme $X$, Pardon\cite{P1} extended the notion of Witt group so that the bilinear forms  take values in certain sheaves of $\oo_X$-modules. In this section, we will review his construction.

Let $X$ be a scheme of dimension $n$, and $\C$ a coherent sheaf of $\oox$-modules with an injective resolution
$$0\ra \C\ra\E^0\rf{d^0}\E^1\rf{d^1}\E^2\ra\cdots\ra\E^n\rf{d^n}0,$$
where $\E^p=\coprod_{x\in X^{(p)}}i\ds E(\kx)$, $E(\kx)$ is the injective hull of the residue class field at $x$ viewed as a constant sheaf on $\bar{x}$, and $i:\bar{x}\hra X$ is the inclusion.
Such $\C$ is called a \itk{canonical sheaf} for $X$ \cite[Chapter 3]{brun}. It is unique up to tensor product with a locally free sheaf of rank 1. Not every scheme admits a canonical sheaf, but every regular scheme does, and $\oo_X$ is a canonical sheaf in such case. Henceforth, unless otherwise stated, we will assume that $X$ is a regular scheme.
Set $\V^p:=\krn d^p$.

\df{$\cm^p(X)$ is the category of CM $\oo_X$-modules of codimension $p$. 

$Q(\cm^p(X);\C)$ is the category of isometry classes of \ns\ symmetric bilinear forms
$$\phi:\M\times\M\ra\V^p,$$
where $\M\in\cm^p(X)$. 
$(\M,\phi)\in Q(\cm^p(X);\C)$ is called a \itk{lagrangian} if there is a submodule $\N\sst\M$ \st $\N,\M/\N\in\cm^p(X)$, $\phi|_{\N\times\N}=0$, and the induced pairing
$$\N\times(\M/\N)\ra\V^p$$
is nonsingular.
}

$Q(\cm^p(X);\C)$ is a semigroup, where the addition defined by the orthogonal sum. The corresponding Grothendieck group modulo the subgroup generated by lagrangians is denoted by $W(\cm^p(X);\C)$.

$\cm^n(X)$ is the category of sheaves of  $\oox$-modules of finite length, and $\cm^0(X)$ is the category of coherent sheaves of locally free $\oox$-modules. Hence, $W(\cm^0(X);\oox)$ is the Witt group $W(X)$ defined by Knebusch \cite{Kn2}. Based on this observation, we will use the notation
$$W(X;\C):=W(\cm^0(X);\C),$$
and if $\oox$ is used as the canonical module, we will suppress it:
$$W(\cm^p(X)):=W(\cm^p(X);\oox)$$

Now we will construct the so-called \itk{lattice map},
$$\Ll^p:\coprod_{x\in X^{(p)}}W(\cm^p(X_x);\C\sx)\dra W(\cm^{p+1}(X);\C).$$

\df{If $x\in X$, let $i_x:\bar{x}\hra X$ be the inclusion map. Suppose that for each $x\in X^{(p)}$, we are given $[M_x,\phi]\in W(X_x;\C\sx)$. Let
$$\phi:=\coprod_{x\in X^{(p)}}i_{x*}\phi\sx,\qqq\N:=\coprod_{x\in X^{(p)}} i_{x*}M_x,$$
where the $\oo_{X,x}$-module $M_x$ is viewed as a constant sheaf on $\bar{x}$. 
An $\oo_X$-submodule $\M\sst \N$ is called a \itk{lattice} if $\M\in\cm^p(X)$ and $\M_x=M\sx\fall x\in X^{(p)}$. The lattice is \itk{integral} \wrt\ $\phi$ if $\phi(\M\times\M)\sst\V^p$. If $\M$ is an integral lattice for $(\N,\phi)$, its \itk{dual lattice} is an $\oox$-module $\M'$ defined for each affine open subset $U\sst X$ by
$$\M'(U)=\setst{n\in \N(U)}{\phi(U)(n,\M(U))\sst\V^p(U)}.$$
}

If $\M$ is an integral lattice for $(\N,\tau)$ and $\M'$ is its dual lattice, there is a \wdf\ bilinear form
$$\bar{\phi}:\frac{\M'}{\M}\times\frac{\M'}{\M}\dra\V^{p+1}$$
given by $\bphi(\bar{m}'_1,\bar{m}'_2)= d^p(\phi(m'_1,m'_2))$ for each affine open subset $U\sst X$. Pardon\cite{bigP} proved that there is a \wdf\ map
$$\Ll^p:\coprod_{x\in X^{(p)}}W(\cm^p(X_x);\C\sx)\dra W(\cm^{p+1}(X);\C),\qq[\N,\phi]\mt[\M'/\M,\bar{\phi}].$$
Unfortunately, the integral lattice does not exist in general. To get around this problem,  Pardon relaxed the condition $\M\in\cm^p(X)$ to a weaker condition $\M\in\Ss^p_2(X)$, where  $\spi(X)$ is the category of coherent sheaves of $\oo_X$-modules $\M$ of codimension $p$ \st
$$\dep_{\oo_{X,x}}\M\sx\geq\min\set{i,\dim_{\oo_{X,x}}\M\sx}\ffall x\in X.$$
He then proved the existence of an $\Ss^p_2$-lattice. There is an inclusion $\cm^p(X)\sst\spi(X)$, which is an equality if $\dm X\leq p+i$. Hence, if $\dm X=2$, then $\cm^p(X)=\spt(X)$, so CM lattice always exists in this case. 

Pardon's original proof \cite[3.12]{bigP} contains an error, which prevents its application to non-affine schemes. Here we include a proof which works for any Gorenstein schemes, provided  that $p=0$ :

\prop{Let $X$ be a Gorenstein scheme, $x\in X^{(p)}$, and $[N,\psi]\in W(\cm^p(X_x);\C_x)$. If $X$ is affine or $p=0$, then there exists an integral lattice for $[N,\psi]$.
}

\pf{Let $i:\bar{x}\hra X$ be the inclusion, and $\M\sst i\ds N$ an $\oox$-submodule \st $\M\sx=N$. Then $\M\in\Ss^p_1(X)$ by \cite[1.19]{bigP}, so the canonical map $\rho:\M\ra\M\ddu$ is injective. On the other hand, since $\M\sx$ is an $\oo_{X,x}$-module of finite length,  $\rho\sx:\M\sx\ra\M\ddu\sx$ is bijective. Thus, we have a map
$$\M\ddu\sx\os{\rho\inv\sx}{\us{\sim}{\ra}}\M\sx=N.$$
Since taking stalks at $x$ is a left adjoint to $i\ds$, we obtain an injective map $\M\ddu\hra i\ds N$. By \cite[1.13]{bigP}, $\M\ddu\in\cm^p(X)$, so $\M\ddu$ is a lattice for $[N,\psi]$.

Now consider the composition
$$\theta:\M\times\M\hra i_*N\times i_*N\rff{i_*\psi}i\ds\E^p\sx\hra\E^p\rf{d^p}\E^{p+1}\ras\bopl_{y\in X^{(p+1)}}i_{y*}\E^{p+1}_y,$$
where $i_y:\bar{y}\hra X$ is the inclusion.
If $X=\spec A$, we can always ``clear out denominators" by multiplying $\M$ by some nonzero element $a\in A$, so that $\theta(a\M\times a\M)=0$, i.e., $i\ds\psi(a\M\times a\M)\sst\V^p$. Then $a\M\sst i\ds N$ is an integral lattice.

This is not always possible, however, if $X$ is not affine:  for example, if $X=\mathbb{P}^1_\cx$, the only regular functions on $X$ are constant functions, so one can't clear out denominators by multiplying by a regular function. We will show that if $p=0$, one can construct a subsheaf $\D\sst\M$ using  a Weil divisor that cancels out the poles appearing in the image of $i\ds\psi$, thus giving $\theta(\D\times\D)=0$.

So assume that $p=0$. Then we may assume that $N=K(X)$, viewed as a constant sheaf on $X$, and that $\psi$ is given by multiplication by some $f\in K(X)$. We may take $\M=\oox$ as our (non-integral) lattice. Our aim is to find a submodule $\D\sst\oox$ \st $\D\in\cm^0(X)$ and $\theta(\D\times\D)=0$.  Let
$$(f):=\sum_{y\in X^{(1)}}n_y(f)y\in\dv(X),$$
and
$$(f)^+:=\sum_{y\in X^{(1)}}n_y^+(f)y,\qqq(f)^+:=\sum_{y\in X^{(1)}}n_y^-(f)y,$$
where
$$n^+_y(f):=\begin{cases}n_y(f) & \ \mb{ if }\  n_y(f)\geq0,\\
          0 & \ \mb{ if }\ n_y(f)<0,
         \end{cases}\qq n^-_y(f):=\begin{cases}0 & \ \mb{ if }\   n_y(f)>0,\\
          n_y(f) & \ \mb{ if }\  n_y(f)\leq0.
         \end{cases}$$
Define an $\oo_X$-module $\D$ for each affine open subset $U\sst X$ by
$$\D(U):=\setst{g\in K(X)}{n_y(g)+n_y^-(f)\geq0 \ \fall y\in U\cap X^{(1)}}.$$
Since 
$$\oo_X(U)=\setst{g\in K(X)}{n_y(g)\geq0 \ \fall y\in U\cap X^{(1)}},$$
$\D$ is a subsheaf of $\oo_X$. Moreover, if $y\in X^{(1)}$ and $\pi_y\sst\oo_{X,y}$ is the maximal ideal, then
$$\D_y=\begin{cases}\oo_{X,y} & \mb{ if } y\not\in\supp(f)^-,\\ \pi_y^{|n_y^-(f)|}\oo_{X,y} & \mb{ if }y\in\supp(f)^-.\end{cases}$$
Hence, if $\eta\in X$ is the generic point, then $\D_\eta=\oo_{X,\eta}=K(X)$, and by construction,
$$\theta_y(\D_y\times\D_y)=0\in\E^1_y\ffall y\in X^{(1)}.$$
Hence, $\theta(\D\times\D)=0$. Finally, since $\D$ is locally free, $\D\in\cm^0(X)$.
}

Given $[\M,\phi]\in W(\cm^p(X);\C)$, we may localize at $x\in X^{(p)}$ to obtain $[\M\sx,\phi\sx]\in W(\cm^p(X\sx);\C\sx)$.
Hence, there is a map
$$\K^p: W(\cm^p(X);\C)\ra\coprod_{x\in X^{(p)}}W(\cm^p(X_x);\C\sx).$$
Pardon\cite[3.9, 3.23]{bigP} showed that the sequence
\beq{0\ra W(\cm^p(X);\C)\rf{\K^p}\coprod_{x\in X^{(p)}}W(\cm^p(X_x);\C\sx)\rf{\Ll^p}W(\cm^{p+1}(X);\C)\label{kl}}
is exact. Setting $d^p:=\K^{p+1}\cc\Ll^p$, we thus obtain a complex
\beq{0\ra\coprod_{x\in X^{(0)}}W(\cm^0(X\sx);\C\sx)\rf{d^0}\cdots\ra\coprod_{x\in X^{(n)}}W(\cm^n(X\sx);\C\sx)\rf{d^n}0.\label{love}}
Moreover, he showed \cite[5.1]{bigP} that if $X$ is the spectrum of a regular local ring which is essentially of finite type over a field of characteristic different from 2, then $\Ll^p$ is surjective and therefore (\ref{love}) is acyclic, with $\krn d^0=W(\Ss^0_1(A);\C)$.
To recover the Gersten-Witt complex from (\ref{love}), he makes use of the \dev\ \cite[2.2]{P3}
\beq{\coprod_{x\in X^{(p)}}W(\kappa(x);\Lambda^pN(\maxi\sx)\ten_{\oo\sx}\C\sx)\ras \coprod_{x\in X^{(p)}}W(\cm^p(X\sx);\C\sx),\label{mon}}
where $N(\maxi\sx):=\Hom_{\oo_x}(\maxi\sx/\maxi\sx^2,\kappa(x))$. One may always choose an \iso
$$\Lambda^pN(\maxi\sx)\ten_{\oox}\C\sx\simeq\kx,$$
giving rise to a \itk{non-canonical} \iso 
$$W(\kappa(x);\Lambda^pN(\maxi\sx)\ten_{\oo\sx}\C\sx)\simeq W(\kappa(x)).$$
Thus we obtain the classical Gersten-Witt complex
$$0\ra\coprod_{x\in X^{(0)}}W(\kappa(x))\rf{d^0}\coprod_{x\in X^{(1)}}W(\kappa(x))\rf{d^1}\cdots\ra\coprod_{x\in X^{(n)}}W(\kappa(x))\rf{d^n}0.$$
The local acyclicity of (\ref{love}) implies that we can sheafify it to obtain a flasque resolution.
Define a \itk{Witt sheaf} for each affine open subset $U\sst X$ by
$$\W(X;\C)(U):=\krn\left(W(\kappa(\eta);\C_\eta)\rf{d^0}\coprod_{x\in U\cap X^{(1)}}W(\cm^1(X\sx);\C\sx)\right),$$
where $\eta\in X$ is the generic point. (\ref{love}) then sheafifies to a flasque resolution of $\W(X;\C)$ :
\begin{multline*}
0\ra i_{\eta*}W(\kappa(\eta);\C_\eta)\rf{d^0}\coprod_{x\in X^{(1)}}i_{x*}W(\kappa(x);N(\maxi\sx)\ten_{\oo_{X,x}}\C\sx)\rf{d^1}\cdots\\
\cdots\ra\coprod_{x\in X^{(n)}}i_{x*}W(\kappa(x);N(\maxi\sx)\ten_{\oo_{X,x}}\C\sx)\rf{d^n}0,
\end{multline*}
where $W(\kx;N(\maxi\sx)\ten_{\oo_{X,x}}\C\sx)$ is viewed as a constant sheaf on $\bar{x}$, and $i_x:\bar{x}\hra X$ is the inclusion \cite[0.11]{bigP}.
By the Purity Theorem\cite{OjP}, the stalk of $\W(X;\C)$ at $x\in X$ is $W(X\sx;\C\sx)$, so $\W(X;\C)(U)=W(U)$.
 
\chapter{Hirzebruch surfaces}
\label{hirze}

This chapter is organized as follows: 

In Section \ref{geo}, we introduce the Hirzebruch surface $H_n$, and discuss its geometry. 

In Section \ref{toric}, we define a new complex of Witt groups on 2-dimensional surfaces supported on codimension $p$ tori.
We define these Witt groups in Sections \ref{witt}. In Section \ref{lattice}, we define the boundary maps of this complex.

In Section \ref{deviss}, we show that the Witt group supported on codimension $p$ tori is \isoc\ to the Witt group of the the tori, thus obtaining a complex of Witt groups of tori. Our claim is that this complex is quasi-\isoc\  to the Gersten-Witt complex of $\hn$. We prove the quasi-\iso\ in the next chapter.

\section{Geometry of Hirzebruch surfaces}
\label{geo}

The Hirzebruch surface $H_n$ is a $\PP^1_\cx$ bundle $\pi:H_n\ra\PP^1_\cx$  obtained by projectivizing the line bundle  $\oo_{\PP^1_\cx}\opl\oo_{\PP^1_\cx}(-n)$. It can be constructed as a toric variety by a fan depicted below \cite[p.~7]{fulton}:
\xym{(-1,n) & (1,0) & \\
& *{} \ar[lu] \ar[r] \ar[d] \ar@{}[ru]|(.3){\sigma_1} \ar@{}[rd]|(.3){\sigma_2} \ar@{}[l]^(.2){\sigma_3} & (1,0) \\
& \ar@{<->}[uu]^(.7){\sigma_4} (-1,0) &
}
To each cone $\sigma_i$ corresponds an affine open subset $\usi$ :
\beq{\xymatrix{U_{\sigma_2} \ar@{}[r]|(.35)= & \spec\cx[x,\by] \ar@{<->}[r] \ar@{<->}[rd] \ar@{<->}[d] & \ar@{<->}[ld]|!{[l];[d]}\hole \ar@{<->}[d] \spec\cx[z,w] \ar@{}[r]|(.65)= & U_{\sigma_3}\\
U_{\sigma_1} \ar@{}[r]|(.35)= & \spec\cx[x,y] \ar@{<->}[r]  & \spec\cx[z,\bw] \ar@{}[r]|(.65)= & U_{\sigma_4}
}\label{coord}}
where $z:=1/x$, $w:=1/x^ny$, $\bw:=1/w$, $\by:=1/y$. 
The axes of $U_{\sigma_i}$ constitute four projective lines in $H_n$, denoted by $X,Y,Z,W$:
\beq{\xymatrix{ &  & \ar@{-}[ddd]|W^(.27){0_{zw}} &  \\
    &  & & \ar@{-}[lll]|Z_(.75){0_{x\by}} \\
      \ar@{-}[rrr]|X_(.75){0_{ z\bw}} &  & &  \\
&\ar@{-}[uuu]|Y^(.27){0_{xy}} & &
}\label{axes}}
$X$ and $Z$ are sections of $\pi:H_n\ra\PP^1_\cx$.
$Y$ and $W$ are fibres over $\PP^1_\cx$, hence have trivial normal bundles, while the normal bundles of $X$ and $Z$ are ``twisted" by the integer parameter $n$. When $n=0$, there is no twist, and $H_0\simeq\PP^1_\cx\times\PP^1_\cx$.

Fern\'{a}ndez-Carmena \cite[3.4]{fer} showed that the Witt group of a smooth complex surface is a birational invariant. Hence,
\beq{W(H_n)=W(\PP^2_\cx)=\Z/2\ffall n\in\Z.\label{fern}}
We will adopt the following notations:
$$H^1_n:=X\cup Y\cup Z\cup W,\qqq H^2_n:=\set{0_{xy},0_{zw},0_{x\by},0_{z\bw}}.$$
$$T_X:=X-\set{0_{xy},0_{z\bw}},\qqq T_Y:=Y-\set{0_{xy},0_{x\by}},$$
$$T_Z:=Z-\set{0_{zw},0_{x\by}},\qqq T_W:=W-\set{0_{zw},0_{z\bw}},$$
$$L:=X\cup Z,\qqq N:=Y\cup W,$$
$$T_L:=T_X\cup T_Z,\qqq T_N:=T_Y\cup T_W.$$
Note that $T_L=H^1_n-N$, $T_N=H^1_n-L$, and $H^2_n=L\cap N$.

\section{A complex of Witt groups supported on tori}
\label{toric}

Let $H=H_n$ for some $n\in\Z$, and $\eta\in H$ the generic point. The Gersten-Witt complex of $H$ is given by
\beq{0\ra W(\kappa(\eta)) \rf{d^0}  \coprod_{x\in H^{(1)}}W(\cm^1(H_x))  \rf{d^1}  \coprod_{x\in H^{(2)}}W(\cm^2(H\sx))\ra0.\label{gerhir}}
In Pardon's construction \cite{bigP},
$$d^0:W(\kappa(\eta))\rf{\Ll^0}W(\cm^1(H))\rf{\K^1}\coprod_{x\in H^{(1)}}W(\cm^1(H_x)),$$
$$d^1:\coprod_{x\in H^{(1)}}W(\cm^1(H_x))\rf{\Ll^1}W(\cm^2(H))\us{\sim}{\rf{\K^2}}\coprod_{x\in H^{(2)}}W(\cm^2(H_x)).$$
We claim (Proposition \ref{tak}, Corollary \ref{ttt}) that there is a quasi-\isoc\ complex
\beq{0\ra W(H-H^1)\rf{d^0_H} W(\cm^1_{T_L}(H-N))\coprod W(\cm^1_{T_N}(H-L)) \rf{d^1_H}W(\cm^2_{H^2}(H))\ra0,\label{toricComplex}}
where $d^0_{H}$ is the compositioin
\beq{W(H-H^1)\rf{\Ll^0_{H^1}} W(\cm^1_{H^1}(H))\rf{\K^1_{H^1}}W(\cm^1_{H^1-N}(H-N))\coprod W(\cm^1_{H^1-L}(H-L)),\label{tor2}}
where $\K^1_{H^1}:=\K^1_N\coprod\K^1_L$,
\begin{eqnarray}
\K^1_N:W(\cm^1_{H^1}(H))&\ra& W(\cm^1_{H^1-N}(H-N)),\\
\K^1_L:W(\cm^1_{H^1}(H))&\ra& W(\cm^1_{H^1-L}(H-L)),
\end{eqnarray}
and $d^1_H:=\Ll^1_{T_L}\coprod\Ll^1_{T_N}$, where
\begin{eqnarray}
\Ll^1_{T_L}:W(\cm^1_{T_L}(H-N))&\ra& W(\cm^2_{L\cap N}(H)),\label{tor3}\\
\Ll^1_{T_N}:W(\cm^1_{T_N}(H-L))&\ra& W(\cm^2_{L\cap N}(H)).\label{tor4}
\end{eqnarray}
Recall that $L\cap N=H^2$. $\K^1_{H^1}$ is an excision map induced by restriction, and it is injective because $\K^1$ is injective \cite[3.9]{bigP}.

In the next couple of sections, we define the Witt groups with support and the lattice maps $\Ll^0_{H^1}, \Ll^1_{T_L}, \Ll^1_{T_N}$.

\section{Witt groups with support}
\label{witt}

If $Y$ is a closed subscheme of $X$, $\cm^p_Y(X)$ will denote the category of  coherent sheaves of CM $\oo_X$-modules of codimension $p$ supported on $Y$. Also, if $\V$ is a sheaf of $\oox$-modules, then $\V_Y$ will denote the sheaf of $\oox$-modules defined for each affine open subset $U\sst X$ by
$$\V_Y(U):=\bigcup_{i=1}^\infty(0:\II(Y)(U)^i)_{\V(U)},$$
where $\II(Y)\sst\oox$ is the ideal sheaf of $Y$.

\subsection{$W(\cm^1_{H^1}(H))$}

Let $Q^1_{H^1}(H)$ be the semigroup of isometry classes of \ns\ symmetric bilinear forms
$$\phi:\M\times \M\ra \V^1_{\ooh},$$
where $\M\in\cm^1_{H^1}(H)$. $(\M,\phi)\in Q^1(H)$ is called a \itk{lagrangian} if there is a submodule $\N\sst\M$ \st $\N,\M/\N\in\cm^1_{H^1}(H)$, $\phi|_{\N\times\N}=0$, and the induced pairing
$$\N\times(\M/\N)\ra\V^1_{\ooh}$$
is \ns. $W(\cm^1_{H^1}(H))$ is the Grothendieck group of $Q^1_{H^1}(H)$ modulo the subgroup generated by lagrangians.
Note that the images of the above bilinear maps lie in $\V^1_{\ooh, H^1}$.

\subsection{$W(\cm^1_{T_L}(H-N))$}

Let $Q^1_{T_L}(H-N)$ be the semigroup of isometry classes of \ns\ symmetric bilinear forms
$$\phi:\M\times \M\ra \V^1_{\oo_{H-N}},$$
where $\M\in\cm^1_{T_L}(H-N)$. $(\M,\phi)\in Q^1_{T_L}(H-N)$ is called a \itk{lagrangian} if there is a submodule $\N\sst\M$ \st $\N,\M/\N\in\cm^1_{T_L}(H-N)$, $\phi|_{\N\times\N}=0$, and the induced pairing
$$\N\times(\M/\N)\ra\V^1_{\oo_{H-N}}$$
is \ns. $W(\cm^1_{T_L}(H-N))$ is the Grothendieck group of $Q^1_{T_L}(H-N)$ modulo the subgroup generated by lagrangians. 
Note that the images of the above bilinear maps lie in $\V^1_{\oo_{H-N}, T_L}$.

\subsection{$W(\cm^1_{T_N}(H-L))$}

Let $Q^1_{T_N}(H-L)$ be the semigroup of isometry classes of \ns\ symmetric bilinear forms
$$\phi:\M\times \M\ra \V^1_{\oo_{H-L}},$$
where $\M\in\cm^1_{T_N}(H-L)$. $(\M,\phi)\in Q^1_{T_N}(H-L)$ is called a \itk{lagrangian} if there is a submodule $\N\sst\M$ \st $\N,\M/\N\in\cm^1_{T_N}(H-L)$, $\phi|_{\N\times\N}=0$, and the induced pairing
$$\N\times(\M/\N)\ra\V^1_{\oo_{H-L}}$$
is \ns. $W(\cm^1_{T_N}(H-L))$ is the Grothendieck group of $Q^1_{T_N}(H-L)$ modulo the subgroup generated by lagrangians. 
Note that the images of the above bilinear maps lie in $\V^1_{\oo_{H-L}, T_N}$.

\subsection{$\cm^2_{H^2}(H)$}

Let $Q^2_{H^2}(H)$ be the  semigroup of isometry classes of \ns\ symmetric bilinear forms
$$\phi:\M\times \M\ra \V^2_{\ooh},$$
where $\M\in\cm^2_{H^2}(H)$. $(\M,\phi)\in Q^2_{H^2}(H)$ is called a \itk{lagrangian} if there is a submodule $\N\sst\M$ \st $\N,\M/\N\in\cm^2_{H^2}(H)$, $\phi|_{\N\times\N}=0$, and the induced pairing
$$\N\times(\M/\N)\ra\V^2_{\ooh}$$
is \ns. $W(\cm^2_{H^2}(H))$ is the Grothendieck group of $Q^2_{H^2}(H)$ modulo the subgroup generated by lagrangians.
Note that the images of the above bilinear maps lie in $\V^2_{\ooh, H^2}$.

\section{Lattice maps with support}
\label{lattice}

We will use $\ooh$ for our canonical sheaf for $H$. The construction is essentially the same as Pardon's \cite{bigP}. All we are doing here is to show that his construction works even with the additional support condition.

\subsection{$\Ll^0_{H^1}:W(H-H^1)\ra W(\cm^1_{H^1}(H))$}

Let $[\N,\psi]\in W(H-H^1)$, so that
$$\psi:\N\times\N\ra\oo_{H-H^1},$$
where $\N\in\cm^0(H-H^1)$. Let $i:H-H^1\hra H$ be the inclusion. An $\ooh$-submodule $\M\sst i\ds\N$ is a \itk{lattice} if $\M\in\cm^0(H)$ and $\M|_{H-H^1}=\N$. The lattice is \itk{integral} \wrt\ $\psi$ if $(i\ds\psi)(\M\times\M)\sst\ooh$. 
\xym{\M\times\M \ar@{}[r]|{\sst} \ar@{-->}[d] & i\ds\N\times i\ds\N \ar[d]^{i\ds\psi}\\
\ooh \ar@{}[r]|{\sst} & i\ds\oo_{H-H^1}
}
If $\M$ is an integral lattice for $[\N,\psi]$, its \itk{dual lattice} is an $\oo_H$-submodule $\M'\sst i\ds\N$ defined for each affine open subset $U\sst H$ by
$$\M'(U)=\setst{n\in i\ds\N(U)}{(i\ds\psi)(U)(n,\M(U))\sst\ooh(U)}.$$
Then there is an induced symmetric bilinear form
$$\bar{\psi}:\frac{\M'}{\M}\times\frac{\M'}{\M}\dra\frac{i\ds\oo_{H-H^1}}{\ooh}$$
given by $\bpsi(\bar{m}'_1,\bar{m}'_2)=d^0((i\ds\psi)(m'_1,m'_2))$ on affine open subsets.

$\M'/\M$ is supported on $H^1$ because $\M'$ locally coincides with $\M$ at every point of $H-H^1$ by the nonsingularity of $\psi$. Moreover, $\M'/\M\in\cm^1(H)$ by \cite[1.2]{P3}.

Before we prove our main propositions (Proposition \ref{main3}, Proposition \ref{main4}), we need some lemmas:

\lemm{In the notation as above, 
$$i\ds\oo_{H-H^1}/\ooh=\V^1_{\ooh}.$$
\label{didy}
}

\pf{$\V^1_{\ooh}=\frac{K(H)}{\ooh}$, where the function field $K(H)$ of $H$ is viewed as a constant sheaf on $H$, and for every affine open subset $U\sst H$, $\V^1_{\ooh,  H^1}(U)\sst (K(H)/\ooh)(U)$ is the subset of sections with poles along $H^1\cap U$:
$$\V^1_{\ooh,H^1}(U)=\bigcup_{j=1}^\infty(0:\II(H^1)^j(U))_{(i\ds K(H)/\ooh)(U)}=\frac{i\ds\oo_{H-H^1}}{\ooh}(U)\sst\frac{K(H)}{\ooh}(U).$$
}

\lemm{Let $[\M,\phi]\in W(\cm^{p+1}_{H^1}(X);\C)$, and $\N\sst \M$ a totally isotropic submodule \st $\N, \N\pp/\N\in\cm^{p+1}(X)$. If the induced bilinear maps
$$\al:\N\times \M/\N\pp\ra\V^{p+1},\qquad \be:\N\pp/\N\times \N\pp/\N\ra\V^{p+1}$$
are \ns\ and $[\N\pp/\N,\be]$ is a lagrangian, then $[\M,\phi]$ is a lagrangian.
\label{oba}
}

\pf{Let $\K/\N\sst \N\pp/\N$ be a sublagrangian, where $\N\sst \K\sst \N\pp$. We will show that $\K\sst \M$ is a sublagrangian. Since $\N,\K/\N\in\cm^{p+1}(X)$, we have $\K\in\cm^{p+1}(X)$ \cite[1.2]{P3}. Being a sublagrangian, $\K/\N\sst \N\pp/\N$ is totally isotropic, so $\K\sst \M$ is totally isotropic, as well. Let
$$\bar{\phi}:\K\times \M/\K\ra\V^{p+1}$$
be the induced pairing. We will show that the induced map
$$\ad\dg\bphi:\M/\K\ra\shom(\K,\V^{p+1})$$
is bijective. This would imply $\M/\K\in\cm^{p+1}(X)$ \cite[1.6a]{P3}, and that $\K\sst\M$ is a sublagrangian, finishing the proof.

To this end, note that there is a short \esq
$$0\ra\frac{\N\pp/\N}{\K/\N}\rf{j}\frac{\M}{\K}\ra\frac{\M}{\N\pp}\ra0.$$
Taking $(-)\du\equiv\shom(-,\V^{p+1})$, we obtain a \cd
\xym{0  & \ar[d]_{\ad\bar{\be}} \K/\N \ar[l] & \ar[d]_{\ad\bar{\phi}} \K \ar[l] & \ar[d]^{\ad\al} \N \ar[l]_i^{\inc} & \ar[l] 0\\
& ((\N\pp/\N)/(\K/\N))\du & (\M/\K)\du \ar[l]_(.35){j\du} & (\M/\N\pp)\du \ar[l] & \ar[l] 0
}
where the rows are exact.  Since $\ad\bar{\be}$ and $\ad\al$ are \isos, $j\du$ is surjective, hence $\ad\bar{\phi}$ is  bijective. Reflexivity of CM modules \cite[1.6a]{P3} then implies that $\ad\dg\bphi$ is also bijective, as desired.
}

\lemm{Let $[\M,\phi]\in W(\cm^{p+1}_{H^1}(X);\C)$, and $\N\sst \M$ a totally isotropic submodule \st $\N, \N\pp/\N\in\cm^{p+1}(X)$. If the induced bilinear maps
$$\al:\N\times \M/\N\pp\ra\V^{p+1},\qquad \be:\N\pp/\N\times \N\pp/\N\ra\V^{p+1}$$
are \ns, then $[\M,\phi]=[\N\pp/\N,\be]$ in $W(\cm^{p+1}_{H^1}(X);\C)$.
\label{oba2}
}

\pf{We will show that $[\M,\phi]-[\N\pp/\N,\be]=[\M\opl(\N\pp/\N),\phi\opl(-\be)]$ is a lagrangian. Let $\N_1\equiv \N\opl0 \sst \M\opl(\N\pp/\N)$. Then $\N_1\in\cm^{p+1}_{H^1}(X)$, and 
$$\N_1\pp=\N\pp\opl(\N\pp/\N)\sst \M\opl(\N\pp/\N),\qquad\frac{\M\opl(\N\pp/\N)}{\N_1\pp}\simeq \M/\N\pp.$$
Since $\N\sst \N\pp$, we have $\N_1\sst\N_1\pp$ and  
$$\N_1\pp/\N_1\simeq(\N\pp/\N)\opl(\N\pp/\N)\in\cm^{p+1}_{H^1}(X).$$
The induced bilinear maps
$$\N_1\times\frac{\M\opl(\N\pp/\N)}{\N_1\pp}\ra\V^{p+1},\qquad\N_1\pp/\N_1\times\N_1\pp/\N_1\ra\V^{p+1},$$
are \isoc\ to
$$\al:\N\times(\M/\N\pp)\ra\V^{p+1},$$
$$\be\opl(-\be):(\N\pp/\N)\opl(\N\pp/\N)\times(\N\pp/\N)\opl(\N\pp/\N)\ra\V^{p+1},$$
\resp. Note that the latter is a hyperbolic form, hence a lagrangian. Thus, applying Lemma \ref{oba} with $[\M\opl(\N\pp/\N),\phi\opl(-\be)]$ in place of $[\M,\phi]$, and $\N_1$ in place of $\N$, it follows that $[\M\opl(\N\pp/\N),\phi\opl(-\be)]$ is a lagrangian.
}

We are now ready to prove one of our main propositions in this chapter:

\prop{There is a \wdf\ map
$$\Ll^0_{H^1}:W(H-H^1)\dra W(\cm^1_{H^1}(H)),\qqq[\N,\psi]\mt[\M'/\M,\bpsi],$$
where $\M$ is an integral lattice for $[\N,\psi]$.
\label{main3}
}

\pf{We have seen above that $\M'/\M\in\cm^1_{H^1}(H)$.
We must show that 1) $\bpsi$ is \ns, 2) $[\M'/\M,\bpsi]$ is independent of the choice of $\M$, and 3) if $[\N,\psi]$ is a lagrangian, so is $[\M'/\M,\bpsi]$.

1) For the nonsingularity of $\bpsi$, we need to show that the induced map
$$\ad\bpsi:\M'/\M\ra\shom(\M'/\M,i\ds\oo_{H-H^1}/\ooh)$$
is bijective.
Injectivity is clear from the definition of dual lattice.
For surjectivity, let $\be\in\shom(\M'/\M,i\ds\oo_{H-H^1}/\ooh)$. The short \esq
$$0\ra\M\ra\M'\ra\M'/\M\ra0$$
induces an \esq\ \cite[1.6b]{P3}
$$0\la (\M'/\M)\hat{~}\lf{\del}\tilde{M}\la\tilde{\M}'\la0.$$
where $(-)\tilde{~}:=\shom(-,\ooh)$, $(-)\hat{~}:=\shom(-,\V^1(H))$.
If $\be=\del(\al)$, \tcd
\xym{0 \ar[r] & \ar[r] \M \ar[d]_\al & \ar[r] \M' \ar[d]_{\al'} & \ar[d]^\be \ar[r] \M'/\M & 0 \\
0 \ar[r] & \ooh \ar[r] & i\ds\oo_{H-H^1} \ar[r] & \V^1_{\ooh, H^1} \ar[r] & 0
}
where the rows are exact (Lemma \ref{didy}).
$\al'$ is given on an affine open subset $U\sst H$ by $\al'(U)=(i\ds\psi)(U)(m,-)$ for some $m\in\M(U)$. 
Hence, 
$$\be(U)=\bpsi(U)(\bar{m}',-)=\ad\bpsi(U)(\bar{m}')$$
for some $\bar{m}'\in(\M'/\M)(U)$, i.e., $\ad\bpsi$ is surjective.

2) Now suppose that $\M_1,\M_2\sst i\ds\N$ are two integral lattices. Then $\M_1\cap\M_2$ is also an integral lattice, so we may assume $\M_1\sst\M_2$. Then there are inclusions
$$\M_1\sst\M_2\sst\M'_2\sst\M'_1\sst i\ds\N.$$
Hence, $\M_2/\M_1\sst \M'_1/\M_1$ is a totally isotropic subspace, and 
$$(\M_2/\M_1)\pp=\M'_2/\M_1\sst \M'_1/\M_1.$$
Hence, there are \isos
$$(\M_2/\M_1)\pp/(\M_2/\M_1)\simeq \M'_2/\M_2,\qquad(\M'_1/\M_1)/(\M_2/\M_1)\pp\simeq \M'_1/\M'_2,$$
and the induced bilinear maps
$$\al:\M_2/\M_1\times(\M'_1/\M_1)/(\M_2/\M_1)\pp\ra\V^1(H),$$
$$\be:(\M_2/\M_1)\pp/(\M_2/\M_1)\times(\M_2/\M_1)\pp/(\M_2/\M_1)\ra\V^1(H),$$
are \isoc\ to the bilinear forms
$$\al':\M_2/\M_1\times \M'_1/\M'_2\ra\V^1(H),$$
$$\be':\M'_2/\M_2\times \M'_2/\M_2\ra\V^1(H),$$
\resp. We have shown the nonsingularity of $\be'$ in the first part of the proof. Similar proof shows that $\ad\al'$ is surjective. The reflexivity of CM modules \cite[1.6]{P3} implies that $\M_1''=\M_1$, which in turn implies injectivity of $\ad\al'$. Hence, $\al'$ and $\be'$ are \ns, and therefore $\al$ and $\be$ are \ns.
Since $\M_2/\M_1,\M'_2/\M_2\in\cm^1(H)$ \cite[1.2]{P3}, the result then follows from applying Lemma \ref{oba2} to $[\M'_1/\M_1,\bpsi]\in W(\cm^1_{H^1}(H))$ and $\M_2/\M_1\sst \M'_1/\M_1$.

3) Let $[\N,\psi]\in W(H-H^1)$ be a lagrangian. We will show that $[\M'/\M,\bpsi]=0\in W(\cm^1_{H^1}(H))$.
Since the localization map
$$W(\cm^1_{H^1}(H))\ra\coprod_{x\in H^{(1)}}W(\cm^1(H\sx))$$
is injective \cite[3.9]{bigP}, it suffices to show that
$$[\M'\sx/\M\sx,\bpsi\sx]=0 \in W(\cm^1(H\sx))\ffall x\in H^{(1)}.$$
Let $\I\sst \N$ be a sublagrangian, and define 
$$\G:=\krn(\M\rightarrowtail i\ds\N\ra i\ds(\N/\I)).$$
Since $\M$ is an integral lattice, there is an induced bilinear map
$$i\ds\psi:\M\times \M\ra\ooh.$$
Since $\G\sst i\ds\I$, the submodule $\G\sst \M$ is totally isotropic \wrt\ $i\ds\psi$, and there is an induced pairing
$$\al:\G\times \M/\G\ra\ooh.$$
Let $(-)\du:=\shom(-,\ooh)$. Since $\M'\simeq\M\du$ \cite[3.16]{bigP}, \tcd
\xym{& 0 \ar[d] & \ar[d] 0 & 0 \ar[d] \\
0 \ar[r] & \ar[r]^j \ar[d]_{\ad\al} \G & \ar[d]_{\ad i\ds \psi} \ar[r]^\pi \M & \ar[d]^{\ad^\dagger\al} \ar[r] \M/\G & 0\\
0 \ar[r] & \ar[r]^{\pi\du} \ar[d] (\M/\G)\du & \ar[d] \ar[r]^{j\du} \M\du & \ar[d] \G\du\\
0 \ar[r] & \ar[d] \ar@{-->}[r]^{\bar{\pi}\du} \cok\ad\al & \ar[d] \ar@{-->}[r]^{\bar{j}\du} \M'/\M & \cok\ad\dg\al \ar[d] \\
& 0 & 0 & 0
}
where the rows and columns are exact.
Since $\M\in\cm^0(H)$ and  $\M/\G\rat i\ds(\N/\I)\in\cm^0(H)$, we have $\G,\M/\G\in\Ss^0_1(H)$ \cite[1.19]{bigP}. Hence, $\G\sx, \M\sx/\G\sx\in\cm^0(H\sx)\fall x\in H^{(1)}$. By \cite[1.6c]{P3}, the second and third rows are locally exact at every point of $H^{(1)}$. Let $\Ss:=\im\bar{\pi}\du$.
It follows from the \cd\ that $\bpsi|_{\Ss\times\Ss}=0$, and the local exactness of the third row implies that $\Ss\sst\M'/\M$ is a sublagrangian at those points (in $H^{(1)}$).
}

Since $\Ll^0_{H^1}$ is defined in the same way as Pardon's lattice map $\Ll^0$, \tcd
\xym{W(\kappa(\eta)) \ar[r]^(.38){d^0} & \coprod_{x\in H^{(1)}}W(\cm^1(H\sx)) \\
W(H-H^1) \ar[r]^(.28){d^0_H} \ar[u]& \ar[u] W(\cm^1_{T_L}(H-N))\coprod W(\cm^1_{T_N}(H-L))
}
where $d^0=\K^1\cc\Ll^0$,  $d^0_H=\K^1_{H^1}\cc\Ll^0_{H^1}$, and the vertical maps are induced by inclusion. Note that the bottom right has only two Witt components supported on $T_L$ and $T_N$, because the image of the lattice map $\Ll^0_{H^1}$ is supported on $H^1$.

\subsection{$\Ll^1_{T_L}\coprod\Ll^1_{T_N}:W(\cm^1_{T_L}(H-N)\coprod W(\cm^1_{T_N}(H-L)\ra W(\cm^2_{H^2}(H))$}

Let $[\N,\psi]\in W(\cm^1_{T_L}(H-N))$, so that
$$\psi:\N\times\N\ra\V^1_{\oo_{H-N}, T_L},$$
where $\N\in\cm^1_{T_L}(H-N)$. Let $i:H-N\hra H$ be the inclusion.
An $\oo_H$-submodule $\M\sst i\ds\N$ is a \itk{lattice} if $\M\in\cm^1_{L}(H)$ and $\M|_{H-N}=\N$. The lattice is \itk{integral} \wrt\ $\psi$ if $(i\ds\psi)(\M\times\M)\sst \V^1_{\ooh}$.
\xym{\M\times\M \ar@{}[r]|{\sst} \ar@{-->}[d] & i\ds\N\times i\ds\N \ar[d]^{i\ds\psi}\\
\V^1_{\ooh} \ar@{}[r]|(.35){\sst} & i\ds\V^1_{\oo_{H-N}, T_L}
}
Since $\M\sst i\ds\N$, the image of the left vertical map necessarily lies in $\V^1_{\ooh, L}$.

If $\M$ is an integral lattice for $[\N,\psi]$, its \itk{dual lattice} is an $\ooh$-submodule $\M'\sst i\ds\N$ defined for each affine open subset $U\sst H$ by
$$\M'(U)=\setst{n\in i\ds\N(U)}{(i\ds\psi)(U)(n,\M(U))\sst \V^1_{\ooh, L}(U)}.$$
Applying Bass's theorem \cite[2.5]{bass} which states that minimal injective resolution is preserved under taking annihilator of a regular element, we obtain an exact sequence
$$0\ra \V^1_{\ooh, L}\ra \E^1_{\ooh, L}\rf{d^1} \E^2_{\ooh,  L}\ra0,$$
which is a minimal injective resolution of $\V^1_{\ooh, L}$ as a sheaf of $\oo^L_H$-modules, where $\oo^L_H$ is the $L$-adic completion of $\ooh$. Thus, there is an induced map
$$\bar{\psi}:\frac{\M'}{\M}\times\frac{\M'}{\M}\dra\V^2_{\ooh}.$$
$\M'$ locally coincides with $\M$ at every point of $H-N$. Also, being submodules of $i\ds\N$, $\M'$ and $\M$ are supported on $L$. Hence, $\M'/\M$ is supported on $N\cap L=H^2$, and the image of $\bpsi$ lies in $\V^2(H)_{H^2}$. Moreover, $\M'/\M\in\cm^2(H)$ by \cite[1.2]{P3}.

Our next main proposition of this chapter can be proved in the same way as Proposition \ref{main3}:

\prop{There are \wdf\ maps
\bea{\Ll^1_{T_L}:W(\cm^1_{T_L}(H-N))&\dra& W(\cm^2_{H^2}(H)),\qq[\N,\psi]\mt[\M'/\M,\bpsi],\\
\Ll^1_{T_N}:W(\cm^1_{T_N}(H-L))&\dra& W(\cm^2_{H^2}(H)),\qq[\N,\psi]\mt[\M'/\M,\bpsi].
}\label{main4}}

Let $d^1_H:=\Ll^1_{T_L}\coprod\Ll^1_{T_N}$, so that
$$d^1_H:W(\cm^1_{T_L}(H-N))\coprod W(\cm^1_{T_N}(H-L))\ra W(\cm^2_{H^2}(H)).$$

Since $\Ll^1_{T_L}$ and $\Ll^1_{T_N}$ are defined in the same way as Pardon's lattice map $\Ll^1$, \tcd
\xym{\coprod_{x\in H^{(1)}}W(\cm^1(H\sx)) \ar[r]^{d^1} & \coprod_{x\in H^{(2)}}W(\cm^2(H\sx))\\
\ar[u] W(\cm^1_{T_L}(H-N))\coprod W(\cm^1_{T_N}(H-L)) \ar[r]^(.67){d^1_H} & \ar[u] W(\cm^2_{H^2}(H))
}
where $d^1=\K^2\cc\Ll^1$.
Note that the Witt group on the bottom right is supported on four points $0_{xy},0_{zw},0_{x\by},0_{z\bw}\in H^{(2)}$.

\section{D\'{e}vissage}
\label{deviss}

Let $X$ be a scheme, and $x\in X^{(p)}$ a point of codimension $p$. In its original form, d\'{e}vissage  \cite[2.2]{P3} states that \tis\ (see (\ref{mon}) for canonical version)
$$W(\kx)\ras W(\cm^p(X\sx)).$$
In order to identify the complex (\ref{toricComplex}) with a complex of Witt groups of tori, we need \isos\ of the form
\beq{W(T_X)\ras W(\cm^1_{T_X}(H_n-N)).\label{devv}}
There is certainly such a map induced by inclusion. However, unlike $\cm^p(X\sx)$, the sheaves in $\cm^1_{T_X}(\hn-N)$ are not of finite length, so \dev\ cannot be applied in its original form.  In this section, we show that the map is still an \iso. First we need a lemma:

\lemm{Let
$$0\ra\E^0\ra\E^1\ra\E^2\ra0$$
be a minimal injective resolution of $\cx[x,y]$. Then
$$\E^1_y:=\bigcup_{i=1}^\infty(0:y^i)_{\E^1}\simeq\frac{\cx(x)[y,1/y]}{\cx(x)[y]}.$$
\label{hew}
}

\pf{Let $\V^p\sst\E^p$ denote the $p$-th cosyzygy. We have $\E^0=\cx(x,y)$, $\V^1=\frac{\cx(x,y)}{\cx[x,y]}$, and
$$\V^1_y:=\bigcup_{i=1}^\infty(0:y^i)_{\V^1}\simeq\frac{\cx[x,y,1/y]}{\cx[x,y]}.$$
Hence, there is an inclusion $\V^1_y\hra\frac{\cx(x)[y,1/y]}{\cx(x)[y]}$. We will show that this is an an essential injective extension over $\cx[x][[y]]$. This is an equivalent way of saying that it is an injective hull over $\cx[x][[y]]$ \cite[2.21]{vamos}. Since $\V^1_y\sst\E^1_y$ is an injective hull over $\cx[x][[y]]$ by Bass \cite[2.5]{bass}, the result follows.

$\frac{\cx(x)[y,1/y]}{\cx(x)[y]}$ is an injective $\cx(x)[[y]]$-module, since it is divisible over PID. This in turn implies that it is injective over $\cx[x][[y]]$.

For essentiality, note that every non-zero element of $\frac{\cx(x)[y,1/y]}{\cx(x)[y]}$ can be represented by a finite sum of the form $\sum_{i\geq1}\frac{f_i}{g_iy^i}$, where $f_i,g_i\in\cx[x]\sst\cx[x][[y]]$. By multiplying this by $\prod_{i\geq}g_i$, one obtains a non-zero element in $\V^1_y=\frac{\cx[x,y,1/y]}{\cx[x,y]}$. This implies that $\V^1_y\hra\frac{\cx(x)[y,1/y]}{\cx(x)[y]}$ is an essential extension.
}

\prop{There are \isos
\bea{W(T_X)\ras W(\cm^1_{T_X}(H_n-N)),&\quad& W(T_Z)\ras W(\cm^1_{T_Z}(H_n-N))\\
W(T_Y)\ras W(\cm^1_{T_Y}(H_n-L)),&\quad& W(T_W)\ras W(\cm^1_{T_W}(H_n-L))
}
induced by inclusion.
\label{learn}
}

\pf{We will only prove the first \iso. The proofs for the other \isos\ are similar.

$T_X$ can be covered by two affine open subsets, $U=\spec\cx[x,1/x]$ and $V=\spec\cx[z,1/z]$, glued together via  $x\leri 1/z$. If $[\M,\phi]\in W(T_X)$, then
\bea{\phi(U):\M(U)\times\M(U)&\ra&\V^0_{\oo_{T_X}}(U)=\cx[x,1/x]\\
\phi(V):\M(V)\times\M(V)&\ra&\V^0_{\oo_{T_X}}(V)=\cx[z,1/z]}
where $\M(U)$ is a free $\cx[x,1/x]$-module, and $\M(V)$ is a free $\cx[z,1/z]$-module. 

On the other hand, viewed as a subset of $H_n-N$, $T_X$ can also be covered by two affine open subsets $U_1=\spec\cx[x,y,1/x]$ and $U_4=\spec\cx[z,\bw,1/z]$, glued together via $x\leri1/z$ and $y\leri z^n\bw$. If $[\N,\psi]\in W(\cm^1_{T_X}(H_n-N))$, then
\bea{\N(U_1)\times \N(U_1)&\rf{\psi(U_1)}&\V^1_{\oo_{H_n-N}, T_X}(U_1)=\bigcup_{i=1}^\infty(0:y^i)_{\frac{\cx(x,y)}{\cx[x,y,1/x]}}=\frac{\cx[x,y,1/x,1/y]}{\cx[x,y,1/x]}\\
\N(U_4)\times \N(U_4)&\rf{\psi(U_4)}&\V^1_{\oo_{H_n-N}, T_X}(U_4)=\bigcup_{i=1}^\infty(0:\bw^i)_{\frac{\cx(z,\bw)}{\cx[z,\bw,1/z]}}=\frac{\cx[z,\bw,1/z,1/\bw]}{\cx[z,\bw,1/z]}
}
where $\N(U_1)$ is a CM $\cx[x,y,1/x]$-module of codimension 1 killed by some power of $y$, and $\N(U_4)$ is a CM $\cx[z,\bw,1/z]$-module of codimension 1 killed by some power of $\bw$.

Let $i:T_X\hra H_n-N$ be the inclusion. There is an injection 
$$j:i\ds\V^0_{\oo_{T_X}}\hra\V^1_{\oo_{H_n-N}, T_X},$$
which is given on the affine charts $U_1$ and $U_4$ by
\bea{j(U_1):\cx[x,1/x]&\os{\frac{1}{y}}{\hra}&\frac{\cx[x,y,1/x,1/y]}{\cx[x,y,1/x]}\\
j(U_4):\cx[z,1/z]&\os{\frac{1}{z^n\bw}}{\hra}&\frac{\cx[z,\bw,1/z,1/y]}{\cx[z,\bw,1/z]}
} 
Thus there is an induced map of Witt groups
$$W(T_X)\dra W(\cm^1_{T_X}(H_n-N)),\qqq[\M,\phi]\mt[i\ds\M,j\cc i\ds\phi].$$
We will show that this map is an \iso. It suffices to show that this is an \iso\ on one of the affine charts, $U_1$.

\Tcd\ of value groups
\xym{& 0\ar[d] & 0 \ar[d] \\
\V^0_{\cx[x,1/x]} \ar@{}[r]|(.55){=} & \cx[x,x\inv] \ar[d] \ar@{^(->}[r]^(.45){1/y} & \frac{\cx[x,y,1/x,1/y]}{\cx[x,y,1/x]} \ar[d] \ar@{}[r]|(.45){=} & \V^1_{\cx[x,y,1/x], y} \\
\E^0_{\cx[x,1/x]} \ar@{}[r]|(.55){=} & \cx(x) \ar[d] \ar@{^(->}[r]^(.45){1/y}  & \ar[d] \frac{\cx(x)[y,1/y]}{\cx(x)[y]} \ar@{}[r]|(.45){=} & \E^1_{\cx(x)[y],y} \\
\V^1_{\cx[x,1/x]} \ar@{}[r]|(.55){=} & \frac{\cx(x)}{\cx[x,1/x]} \ar[d] \ar@{^(->}[r]^(.35){1/y}  & \ar[d] \frac{\cx(x)[y,1/y]}{\cx[x,y,1/x,1/y]+\cx(x)[y]} \ar@{}[r]|(.6){=} & \V^2_{\cx[x,y,1/x], y}  \\
& 0 & 0
}
where the columns are exact. Note that we can make the identification
$$\E^1_{\cx(x)[y],y}=\frac{\cx(x)[y,1/y]}{\cx(x)[y]}$$
by Lemma \ref{hew}. Thus we have an induced \cd\ of Witt groups
\xym{0\ar[d] & 0\ar[d] \\
W(\cm^0(\cx[x,1/x])\ar[d]_{\K^0} \ar[r]^(.45){1/y} & \ar[d]_{\K^1_y}  W(\cm_y^1(\cx[x,y,1/x])) \\
W(\cm^0(\cx(x))) \ar[d]_{\Ll^0} \ar[r]^(.45){1/y} & W(\cm_y^1(\cx(x)[y])) \ar[d]_{\Ll^1_y}  \\
\ar[d] \ar[r];[rd] W(\cm^1(\cx[x,1/x])) \ar[r]^(.45){1/y} &  W(\cm_y^2(\cx[x,y,1/x])) \\
0 & 0
}
Pardon\cite{bigP} showed that the first column is exact. His proof easily extends to the second column (shown below), so the second column is also exact. We will show that the second and third horizontal maps are \isos, which implies that the first horizontal map is also an \iso.

To this end, note that the modules in $\cm^1_y(\cx(x)[y])$ and $\cm^2_y(\cx[x,y,1/x])$ are of finite length. Hence, we may apply \dev\ \cite[2.2]{P3} to reduce the powers of $y$ which annihilate the modules down to 1, thereby representing the forms by the images of the horizontal maps. This proves surjectivity. Injectivity follows from the fact that \dev\ preserves lagrangians.

To prove exactness of the second column, first note that surjectivity of $\Ll^1_y$ follows from surjectivity of $\Ll^0$ and bijectivity of the bottom horizontal map. 
Secondly, injectivity of $\K^1_y$  follows from injectivity of $\K^1$ \cite[3.9]{bigP}. $\Ll^1_y\cc\K^1_y=0$ is clear, and the only non-trivial part is to show that $\krn\Ll^1_y\sst\im\K^1_y$.
Let $[N,\psi]\in \krn\Ll^1_y$, so that
\beq{[M'/M,\bar{\psi}]+[L_1,\phi_1]= [L_2,\phi_2],\label{hyp}}
where $M$ is an integral lattice for $[N,\psi]$, $M'$ is its dual lattice, and $[L_1,\phi_1]$, $[L_2,\phi_2]$ are lagrangians.

First suppose that $[L_1,\phi_1]=0$, so that $[M'/M,\bpsi]$ is a lagrangian. Let $ \bar{K}\sst M'/M$ be a sublagrangian, and $ K\sst M'$ its pullback under the quotient map $ M'\sur M'/M$. Since $\bar{K}$ is a sublagrangian, $\bar{K}= \bar{K}^\perp$ \cite[p.~134]{Kn2}, therefore $K= K'$. The \iso
$$K' \simeq \Hom(K,\V^1_{\cx[x,y,1/x]}) $$
then implies that $K\in\cm^1_y(\cx[x,y,1/x])$ \cite[1.13]{bigP}, and that
$$\psi|_{K\times K}: K\times K\ra\V^1_{\cx[x,y,1/x]}$$
is \ns. Hence, $[ K,\psi|_{K\times K}]\in W(\cm_y^1(\cx[x,y,1/x]))$, 
and clearly 
$$ \K^1_y([K,\psi_{ K\times K}])=[N,\psi],$$
proving $\krn\Ll^1_y\sst\im\K^1_y$.

Next, we show that $[L_1,\phi_1]\in\im \K^1_y$, justifying our assumption. By adding $[L_1,-\phi_1]$ to both sides of (\ref{hyp}), we may assume that $[L_1,\phi_1]$ is a hyperbolic form  (\ref{useful}), so that
$$ L_1= T\opl \tilde{T},$$
where $T\in\cm^2_y(\cx[x,y,1/x])$ and $ \tilde{T}:=\Hom( T,\V^2_{\cx[x,y,1/x]}).$
By \cite[1.6b]{P3}, there exists an \esq
\beq{0\ra   I\ra J\ra  T\ra0,\label{idk1}}
where $ J\in\cm^1_y(\cx[x,y,1/x])$, and a dual exact sequence
\beq{0\ra  J\du\ra  I\du \ra \tilde{T}\ra0,\label{idk2}}
where $J\du:=\Hom(J,\V^1_{\cx[x,y,1/x]})$, $ I\du:=\Hom(  I,\V^1_{\cx[x,y,1/x]})$.
Combining (\ref{idk1}) and (\ref{idk2}) gives an \esq
$$0\ra I\opl J\du\ra J\opl I\du\ra L_1\ra0.$$
Let $ S:=  I\opl J\du$, and
$$\sigma: S\times S\ra \V^1_{\cx[x,y,1/x]}$$
the induced pairing with $\sigma_{|  I\times  I}=\sigma_{|J\du\times J\du}=0$.
Then $[S,\sigma]\in W(\cm^1_y(\cx[x,y,1/x]))$ \cite[1.6a]{P3}, and $\K^1_y([S,\sigma])=[ L_1,\phi_1]$, as desired.
}

\cor{There are canonical \isos
\bea{W(T_X)\coprod W(T_Z) &\ras& W(\cm^1_{T_L}(H_n-N)) \\
W(T_Y)\coprod W(T_W) &\ras& W(\cm^1_{T_N}(H_n-L))
} induced by inclusion.
\label{ttt}
}

\pf{
It is easy to see that if $\M\in\cm^1_{T_L}(H_n-N)$, then $\M=\M_X\opl\M_Z$, where $\M_X:=\cup_{i=1}^\infty(0:\II(X)^i)_\M$, $\M_Z:=\cup_{i=1}^\infty(0:\II(Z)^i)_\M$, and $\II(X),\II(Z)\sst\oo_{H_n}$ are the ideal sheaves of $X$ and $Z$, \resp. Hence, \tis
\bea{W(\cm^1_{T_L}(H_n-N))&\simeq& W(\cm^1_{T_X}(H_n-N))\coprod W(\cm^1_{T_Z}(H_n-N)),\\
  {[} \M,\phi] &\mt& ([\M_X,\phi|_{\M_X}],[\M_Z,\phi|_{\M_Z}])\\
{[}\M,\phi]+[\N,\psi] &\mf&([\M,\phi],[\N,\psi])
}and the result follows from Proposition \ref{learn}.
}

\chapter{Quasi-isomorphism via toric decomposition}
\label{quasi}

In this chapter, we will prove that the toric complex  
\beq{0\ra W(H_n-H_n^1)\rf{d^0_{\hn}} W(H_n^1-H_n^2)\rf{d^1_{\hn}} W(H_n^2)\ra0}
is quasi-\isoc\ to the Gersten-Witt complex  of $\hn$ 
\beq{0\ra \coprod_{x\in H_n^{(0)}}W(\kappa(x))\rf{d^0}\coprod_{x\in H_n^{(1)}}W(\kappa(x))\rf{d^1}\coprod_{x\in H_n^{(2)}}W(\kappa(x))\rf{d^2}0.\label{mon2}}
By Corollary \ref{ttt}, this implies that the complex 
\beq{0\ra W(H-H^1)\rf{d^0_H} W(\cm^1_{T_L}(H-N))\coprod W(\cm^1_{T_N}(H-L)) \rf{d^1_H}W(\cm^2_{H^2}(H))\ra0\label{mon3}}
that we constructed in Chapter \ref{hirze} is quasi-\isoc\ to (\ref{mon2}).

More generally, let $X$ be a toric variety of dimension $n$, where
$$X=X^0\supset X^1\supset\cdots\supset X^n\supset X^{n+1}=\ets$$
is a chain of closures of orbits of the torus action, and $Y^p:=X^p-X^{p+1}$ is a finite disjoint union of $(n-p)$-tori.
Let
$$R^p(X):=\bopl_{x\in X^{(p)}}W(\kappa(x))$$
be the $p$-th term of the Gersten-Witt complex of $X$.

In order to prove the quasi-\iso, we need the following proposition, due to Pardon:

\prop{Let $k$ be a field with $\chr k\neq2$. Then the Gersten-Witt complex of $\Aa^n_k$ is acyclic, and $H^0(\Aa^n_k)=W(\Aa^n_k)$.
\label{torus}
}

\pf{We will prove by induction on $n$. It is trivially true if $n=0$, so assume that $n\geq1$, and that it is true for $n-1$.

If $\pr\in\spec k[x_1]$, denote its fibre under the projection
$$\spec k[x_1,\ldots,x_n]\sur\spec k[x_1]$$
by $F_\pr$. If $\pr\in\spec k[x_1]$ is not the generic point, then $\pr=(f)$, where $f\in k[x_1]$ is an irreducible polynomial. Hence,
$$F\spr=\spec(k[x_1,\ldots,x_n]/(f))=\spec((k[x_1]/(f))[x_2,\ldots,x_n]),$$
and $k[x_1]/(f)$ is an algebraic field extension over $k$.
If $\eta\in\spec k[x_1]$ is the generic point, then $F_\eta=\spec k(x_1)[x_2,\ldots,x_n]$. Note that $\spec k[x_1,\ldots,x_n]$ and $F_\eta$ have the same function field, $k(x_1,\ldots,x_n)$.

Let $A:=\spec k[x_1,\ldots,x_n]$ and $A_1:=\spec k[x_1]$. \Tcd
\xym{0 \ar[r] & \ar[r] \ar[d] W(k[x_1]) & \ar[r] \ar[d] W(k(x_1)) & \ar[r] \ar[d] \coprod_{\pr\in A_1^{(1)}} W(k[x_1]/\pr) & 0 \\
0 \ar[r] & W(A) \ar[r] & W(F_\eta) \ar[r] & \coprod_{\pr\in A_1^{(1)}}W(F\spr) \ar[r] & 0
}
where the vertical maps are induced by inclusion. By Karoubi \cite{kar}, the vertical maps are \isos, and by Pardon \cite{P3}, the first row is exact. Hence, the second row is also exact.

Now, there is a short exact sequence of Gersten-Witt complexes
$$0\ra \coprod_{\pr\in A_1^{(1)}}R\sbu(F\spr)[-1]\ra R\sbu(A)\ra R\sbu(F_\eta)\ra0,$$
where $[-1]$ indicates a degree shift by -1. \Tcd
\xym{&&& 0 \ar[d] \\
& 0 \ar[d] &  & \ar[d]  \coprod_{\pr\in A^{(1)}_1} W(K(F\spr)) \ar[r] & \dots\\
0\ar[r] & \ar[d] W(A) \ar[r]^(.45){\ep_A} & \ar@{}[d]|{\rotatebox{90}{$=$}} W(K(A)) \ar[r]^(.39){d^0_A} & \ar[d] \coprod_{x\in A^{(1)}} W(\kx)  \ar[r] & \dots\\
0\ar[r] & \ar[d] \ar[r]^(.45){\ep_{F_\eta}} W(F_\eta) & W(K(F_\eta)) \ar[r]^(.39){d^0_{F_\eta}}  & \coprod_{x\in F_\eta^{(1)}} W(\kx) \ar[d]  \ar[r] & \dots\\
&  \coprod_{\pr\in A^{(1)}_1} W(F\spr) \ar[d] && 0 \\
& 0
} where columns are exact, and $K(-)$ denotes the function field. Diagram chasing shows that there is an induced map
$$\coprod_{\pr\in A^{(1)}_1}W(F\spr)\dra \coprod_{\pr\in A^{(1)}_1}W(K(F\spr)).$$
By the induction hypothesis, this map is injective, and the first and third rows in the above diagram are exact. Hence, the second row also is exact.
}

The same proof works with Laurent polynomials:

\prop{Let $k$ be a field with $\chr k\neq2$, and $T$ a torus (of any dimension) over $k$. Then the Gersten-Witt complex of $T$ is acyclic, and $H^0(T)=W(T)$.
\label{torus2}
}

We now prove the main proposition of this chapter (see Takeda \cite{takeda} for $K$-theoretic analogue):

\prop{ The complex $W(Y^\bullet)$ is quasi-\isoc\ to $R^{\bullet}(X)$.
\label{tak}
}

\pf{The inclusion $X^{p+1}\hra X^p$ induces a short \esq
\beq{0\ra R^{\bullet}(X^{p+1})[-1]\ra R^{\bullet}(X^p)\ra R^{\bullet}(Y^p)\ra0,\label{keith}}
where $[-1]$ indicates a degree shift by -1. 
Since $Y^p=\coprod_{i}T^p_i$, $R^\bullet(Y^p)$ is acyclic by Proposition \ref{torus2}. Hence, the short \esq\ (\ref{keith}) induces an exact sequence
\beq{0\ra H^0(R^{\bullet}(X^p))\rf{\K^p} H^0(R^{\bullet}(Y^p))\os{\del^p}{\dra} H^0(R^{\bullet}(X^{p+1}))\rf{\J^p} H^1(R^{\bullet}(X^p))\ra0,\label{globe}}
and \isos
$$H^{k-1}(R^{\bullet}(X^{p+1}))\ras H^k(R^{\bullet}(X^p))\ffall k\geq2.$$
The latter gives rise to a chain of \isos
$$H^1(\rb(X^{p-1}))\ras H^2(\rb(X^{p-2}))\ras\cdots\ras H^p(\rb(X^0)),$$
which are induced by the inclusion $\rb(X^{p-i})[-1]\hra\rb(X^{p-i-1})$.

Let $\pl^p:=\K^{p+1}\cc\del^p$. The exactness of the sequence (\ref{globe}) implies that there is a complex
\beq{\cdots\ra H^0(\rb(Y^{p-1}))\rf{\pl^{p-1}}H^0(\rb(Y^p))\rf{\pl^p}H^0(\rb(Y^{p+1}))\ra\cdots.\label{palin}}
By Proposition \ref{torus},
$$H^0(R^\bullet(Y^p))=W(Y^p)=\coprod_iW(T^p_i).$$
Hence, (\ref{palin}) gives a complex $W(Y^\bullet)$.

\Tcd
\xym{& H^0(\rb(Y^{p-1})) \ar[d]_{\del^{p-1}} \ar[rd]^{\pl^{p-1}} & & 0 \ar[d] \\
0 \ar[r] & \ar[d]_{\J^{p-1}} H^0(\rb(X^p)) \ar[r]^{\K^p} & H^0(\rb(Y^p)) \ar[r]^(.47){\del^p} \ar[rd]_(.45){\pl^p} & H^0(\rb(X^{p+1})) \ar[d]^{\K^{p+1}} \\
& \ar[d] H^1(\rb(X^{p-1})) && H^0(\rb(Y^{p+1}))\\
& 0
} where the rows and columns are exact.
Hence, \tcd
\xym{H^0(\rb(X^p)) \ar[r]^(.6){\K^p}_(.6)\sim & \krn \pl^p \\
\ar@{}[u]|{\rotatebox{90}{$\sst$}} \im\del^{p-1} \ar[r]^\sim & \im \pl^{p-1} \ar@{}[u]|{\rotatebox{90}{$\sst$}} 
}
which gives rise to \isos
$$H^p(\rb(X^0))\las H^1(\rb(X^{p-1}))\os{\J^{p-1}}{\us{\sim}{\la}}\frac{H^0(\rb(X^p))}{\im \del^{p-1}}\os{\K^p}{\us{\sim}{\ra}}\frac{\krn \pl^p}{\im \pl^{p-1}}=H^p(W(Y^\bullet)),$$
where the first two \isos\ are induced by inclusion.

Now we will show that this \iso\ is induced by a chain map 
$$\rb(X^0)\dla W(Y^\bullet).$$
By the short exact sequence (\ref{keith}), \tcd
\xym{& & 0 \ar[d] & 0 \ar[d] \\
& 0 \ar[d] & H^0(\rb(X^p))  \ar[d]^{\ep(X^p)} \ar[d] \ar@{-->}[r] &  \ar[d] H^0(\rb(Y^p))  \ar@{-->}@/_3pc/[lld]_(.7){\del^p} \ar[d]^{\ep(Y^p)} \\
& H^0(\rb(X^{p+1})) \ar[d]^{\ep(X^{p+1})} \ar[d] & \ar[r]^\sim R^0(X^p)  \ar[d]^{d^0(X^p)} & R^0(Y^p) \ar[d]^{d^0(Y^p)} \\
0 \ar[r] & R^0(X^{p+1}) \ar[d]^{d^0(X^{p+1})} \ar[r] & R^1(X^p) \ar[d]^{d^1(X^{p})} \ar[r] & R^1(Y^p) \ar[d]^{d^1(Y^{p})} \ar[r] & 0 \\
0 \ar[r] & R^1(X^{p+1}) \ar[d]^{d^1(X^{p+1})} \ar[r] & R^2(X^p) \ar[d]^{d^2(X^{p})} \ar[r] & R^2(Y^p) \ar[d]^{d^2(Y^{p})} \ar[r] & 0 \\
0 \ar[r] & \ar[d] R^2(X^{p+1}) \ar[r] & \ar[d] R^3(X^p) \ar[r] & \ar[d]  R^3(Y^p) \ar[r] & 0 \\
&\vdots&\vdots&\vdots
} where the rows and columns are exact.
By the \iso\ $R^0(X^p)\ras R^0(Y^p)$, there is a map $\lam^p:H^0(R^\bullet(Y^p))\dra R^0(X^p)$, and a \cd
\xym{H^0(\rb(Y^{p-1})) \ar[d]_{d^{p-1}}   \ar@{^(->}[r]^(.55){\lam^{p-1}}  &  R^0(X^{p-1}) \ar[rd]^{d^0(X^{p-1})} \\
H^0(\rb(Y^{p})) \ar[d]_{d^{p}}   \ar@{^(->}[r]^(.55){\lam^{p}} &   R^0(X^{p}) \ar@{^(->}[r] \ar[rd]^(.6){d^0(X^p)} & R^1(X^{p-1})  \ar[rd]^(.6){d^1(X^{p-1})}  \\
H^0(\rb(Y^{p+1}))  \ar@{^(->}[r]^(.55){\lam^{p+1}} &  R^0(X^{p+1}) \ar@{^(->}[r] & R^1(X^p) \ar@{^(->}[r] & R^{2}(X^{p-1})
}
Hence, there are inclusions of chains
$$W(Y^\bullet)\os{\lam^\bullet}{\hra}\rb(X^{p-1})[-p+1]\hra\rb(X^{p-2})[-p+2]\hra\cdots\hra\rb(X^0),$$
which induces the \iso\ $H^p(W(Y^\bullet))\ras H^p(\rb(X^0))$.
}

\chapter{Computations}
\label{compu}

In this chapter, we will compute the boundary maps of the toric complex (\ref{toricComplex}). For the sake of simplicity (and without loss of generality), we will assume that $n\geq0$.

From our choice of affine coordinates (\ref{coord}), we have $H_n-H^1_n=\spec\cx[x,y,1/x,1/y]$.
Then $W(H_n-H^1_n)$ is a $\Z/2$-vector space of dimension 4, generated by the unary forms $\ang{1},\ang{x},\ang{y},\ang{xy}$ \cite[3.11]{kar}. On the other hand, by Corollary \ref{ttt},
\beq{W(\cm^1_{T_L}(H_n-N))\coprod W(\cm^1_{T_N}(H_n-L))\label{horn}}
is generated by 8 basis elements corresponding to the basis elements of 
$$W(T_X),\quad W(T_Y),\quad W(T_Z),\quad W(T_W),$$
each of which is generated by two basis elements \cite[3.9]{kar}. Hence, $d^0_{H_n}$ can be represented by an 8-by-4 matrix. On the other hand, $H_n^2$ consists of four points (refer to the picture (\ref{axes})), and since the Witt group of a point is $W(\cx)=\Z/2$, $W(\cm^2_{H^2_n}(H_n))$ is a $\Z/2$-\vsp\ of dimension 4, so $d^1_{H_n}$ can be represented by a 4-by-8 matrix.

\prop{With above choice of basis, the matrix representation for $d^0_{H_n}$ is given by
$$d^0_{H_{\mb{\scriptsize even}}}=
  \begin{blockarray}{ccccc}
    & \ang{1} & \ang{x} & \ang{y} & \ang{xy} \\
    \begin{block}{c(cccc)}
    \ang{1_x} & 0 & 0 & 1 & 0\\
    \ang{x} & 0 & 0 & 0 & 1 \\
    \ang{1_y} & 0 & 1 & 0 & 0 \\
    \ang{y} & 0 & 0 & 0 & 1 \\
    \ang{1_z} & 0 & 0 & 1 & 0 \\
    \ang{z} & 0 & 0 & 0 & 1 \\
    \ang{1_w} & 0 & 1 & 0 & 0 \\
    \ang{w} & 0 & 0 & 0 & 1 \\
    \end{block}
      \end{blockarray} \quad,\quad
  d^0_{H_{\mb{\scriptsize odd}}}=
  \begin{blockarray}{ccccc}
    & \ang{1} & \ang{x} & \ang{y} & \ang{xy} \\
    \begin{block}{c(cccc)}
    \ang{1_x} & 0 & 0 & 1 & 0\\
    \ang{x} & 0 & 0 & 0 & 1 \\
    \ang{1_y} & 0 & 1 & 0 & 0 \\
    \ang{y} & 0 & 0 & 0 & 1 \\
    \ang{1_z} & 0 & 0 & 0 & 1 \\
    \ang{z} & 0 & 0 & 1 & 0 \\
    \ang{1_w} & 0 & 1 & 0 & 0 \\
    \ang{w} & 0 & 0 & 1 & 0 \\
    \end{block}
      \end{blockarray}\quad.
$$\label{main1}
}

\pf{Let us first determine the entries in the first column. To do this, we need to see where the form $(\oohm,\ang{1})$ is sent to by the composition $\K^1_{H^1_n}\cc\Ll^0_{H^1_n}$. To apply the lattice map $\Ll^0_{H^1_n}$, we need to find an integral lattice for $(\oohm,\ang{1})$. We claim that $\oohn$ is an integral lattice for $(\oohm,\ang{1})$. To see this, we check that the image of the bilinear form
$$\oohn(\usi)\times\oohn(\usi)\rf{\ang{1}}(j\ds\oohm)(\usi)$$
lies in $\V^0_{\oo_{H_n}}(U_{\sigma_i})$ for every $i$, where $j:\hn-\hn^1\hra\hn$ is the inclusion  (see (\ref{coord}) for  the definition of $\usi$.)  For example, on $U_{\sigma_1}$, the form $(\oo_{H_n-H^1_n},\ang{1})$ is given by
$$\cx[x,y,1/x,1/y]\times\cx[x,y,1/x,1/y]\rf{1}\cx[x,y,1/x,1/y],$$
and the image of the bilinear form
$$\oo_{H_n}(U_{\sigma_1})\times\oo_{H_n}(U_{\sigma_1})=\cx[x,y]\times\cx[x,y]\rf{1}\cx[x,y,1/x,1/y]$$
lies in $\V^0_{\oo_{H_n}}(U_{\sigma_1})=\oo_{H_n}(U_{\sigma_1})=\cx[x,y]$. 

To find its dual lattice, note that
$$\oo'_{H_n}(U_{\sigma_1}):=\setst{f\in\cx[x,y,1/x,1/y]}{f\cdot\cx[x,y]\sst\cx[x,y]}=\cx[x,y]=\oohn(U_{\sigma_1}).$$
We can similarly check that $\oo'_{H_n}(\usi)=\oohn(\usi)$ for every $i$, so $\oohn$ is self-dual, resulting in $\Ll^0_{H^1_n}(\ang{1})=0$. Hence, $d^0_{H_n}(\ang{1})=\K^1_{H^1_n}\cc\Ll^0_{H^1_n}(\ang{1})=0$, and the entries in the first column are all 0.

Now let us determine the entries in the second column.

The form $(\oohm,\ang{x})$ is given on the affine charts by
\bea{U_{\sigma_1} \quad :\quad \cx[x,y,1/x,1/y]\times\cx[x,y,1/x,1/y]&\rf{x}&\cx[x,y,1/x,1/y]\\
U_{\sigma_2} \quad :\quad\cx[x,\by,1/x,1/\by]\times\cx[x,\by,1/x,1/\by]&\rf{x}&\cx[x,\by,1/x,1/\by]\\
U_{\sigma_3} \quad :\quad\cx[z,w,1/z,1/w]\times\cx[z,w,,1/z,1/w]&\rf{1/z}&\cx[z,w,1/z,1/w]\\
U_{\sigma_4} \quad :\quad\cx[z,\bw,1/z,1/\bw]\times\cx[z,\bw,1/z,1/\bw]&\rf{1/z}&\cx[z,\bw,1/z,1/\bw]
}
This time, $\oohn$ is not an integral lattice because the image of the bilinear form
$$\oo_{H_n}(U_{\sigma_3})\times\oo_{H_n}(U_{\sigma_3})=\cx[z,w]\times\cx[z,w]\rf{1/z}\cx[z,w,1/z,1/w]$$
does not lie in $\V^0_{\oohn}(U_{\sigma_3})=\ooh(U_{\sigma_3})=\cx[z,w]$. On the other hand,
$$\oo_{H_n}(-W)(U_{\sigma_3})=z\cdot\cx[z,w],$$
and the image of the bilinear form
$$z\cdot\cx[z,w]\times z\cdot\cx[z,w]\rf{1/z}\cx[z,w,1/z,1/w]$$
does lie in $\cx[z,w]$. We can similarly check that the image of the bilinar form
$$\oohn(-W)(\usi)\times\oohn(-W)(\usi)\rf{\ang{x}} (j\ds\oohm)(\usi)$$
lies in $\V^0_{\oo_{H_n}}(\usi)$ for every $i$. Hence, $\oohn(-W)$ is an integral lattice for $(\oohm,\ang{x})$. To find its dual lattice, note that
\bea{\oo_{H_n}(-W)'(U_{\sigma_1})&:=&\setst{f\in\cx[x,y,1/x,1/y]}{f\cdot x\cdot\cx[x,y]\sst \cx[x,y]}\\
&=&\frac{1}{x}\cdot\cx[x,y]\\
&=&\oohn(Y)(U_{\sigma_1}),
}
We can similarly check on the other affine open subsets to conclude that $\oohn(-W)'=\oohn(Y)$.
Hence, $\Ll^0_{H^1_n}(\ang{x})$ is given by
\beq{\frac{\oo_{H_n}(Y)}{\oo_{H_n}(-W)}\times\frac{\oo_{H_n}(Y)}{\oo_{H_n}(-W)}\rf{\ang{x}}\frac{i_*\oo_{H_n-H^1_n}}{\oo_{H_n}}.\label{hel}}
Now we apply the map
$$\K^1_{H^1_n}:W(\cm^1_{H^1_n}(H_n))\ra W(\cm^1_{H^1_n-N}(H_n-N))\coprod W(\cm^1_{H^1_n-L}(H_n-L))$$
by restricting the domain from $H_n$ to $H_n-N$ and $H_n-L$.   (Recall $L:=X\cup Z$, $N:=Y\cup W$.)
Restricting (\ref{hel}) to $H_n-N$ gives zero because
$$\left.\frac{\oohn(Y)}{\oohn(-W)}\right|_{H_n-N}=\frac{\oo_{H_n-N}}{\oo_{H_n-N}}=0.$$
By Corollary \ref{ttt}, this implies that $\K^1_{H^1_n}\cc\Ll^0_{H^1_n}(\ang{x})$ has no component in the subspace generated by $\ang{1_x}, \ang{x}, \ang{1_z}, \ang{z}$, hence the corresponding rows in the second column are 0.

On the other hand, restricting (\ref{hel}) to $H_n-L$ gives
\beq{\frac{\oo_{H_n-L}(T_Y)}{\oo_{H_n-L}(-T_W)}\times\frac{\oo_{H_n-L}(T_Y)}{\oo_{H_n-L}(-T_W)}\rf{\ang{x}}\frac{i_*\oo_{H_n-H^1_n}}{\oo_{H_n-L}}.\label{zz}}
On $(H_n-L)\cap U_{\sigma_1}=\spec\cx[x,y,1/y]$, (\ref{zz}) gives a \cd
\beq{\xymatrix{\frac{\frac{1}{x}\cdot\cx[x,y,1/y]}{\cx[x,y,1/y]}\times\frac{\frac{1}{x}\cdot\cx[x,y,1/y]}{\cx[x,y,1/y]} \ar[r]^(.63){x} & \frac{\cx[x,y,1/x,1/y]}{\cx[x,y,1/y]}\\
\dssimu^{1/x\times1/x} \cx[y,1/y]\times\cx[y,1/y] \ar[r]^(.61)1 & \cx[y,1/y] \ar@{^(->}[u]_{1/x}
}\label{sun2}}
while on $(H_n-L)\cap U_{\sigma_3}=\spec\cx[z,w,1/w]$, it gives a \cd
\xym{\frac{\cx[z,w,1/w]}{z\cdot\cx[z,w,1/w]}\times\frac{\cx[z,w,1/w]}{z\cdot\cx[z,w,1/w]} \ar[r]^(.63){1/z} & \frac{\cx[z,w,1/z,1/w]}{\cx[z,w,1/w]}\\
\dssimu^{1\times1} \cx[w,1/w]\times\cx[w,1/w] \ar[r]^(.61)1 & \cx[w,1/w] \ar@{^(->}[u]_{1/z}
}
We therefore conclude that $\K^1_L\cc\Ll^0_{H_n}(\ang{x})=\ang{1_y}+\ang{1_w}$. Hence, in the second column, the rows corresponding to $\ang{1_y}$ and $\ang{1_w}$ are 1, and the rows corresponding to $\ang{y}$ nad $\ang{w}$ are 0.

So far our results didn't depend on $n$. Now we will see that the third and fourth columns do depend on $n$.

To determine the third column, consider the form $(\oohm,\ang{y})$, given on the affine charts by
\bea{U_{\sigma_1} \quad :\quad \cx[x,y,1/x,1/y]\times\cx[x,y,1/x,1/y]&\rf{y}&\cx[x,y,1/x,1/y]\\
U_{\sigma_2} \quad :\quad\cx[x,\by,1/x,1/\by]\times\cx[x,\by,1/x,1/\by]&\rf{1/\by}&\cx[x,\by,1/x,1/\by]\\
U_{\sigma_3} \quad :\quad\cx[z,w,1/z,1/w]\times\cx[z,w,,1/z,1/w]&\rf{z^n/w}&\cx[z,w,1/z,1/w]\\
U_{\sigma_4} \quad :\quad\cx[z,\bw,1/z,1/\bw]\times\cx[z,\bw,1/z,1/\bw]&\rf{z^n\bw}&\cx[z,\bw,1/z,1/\bw]
}
We use the same argument as above to conclude that $\oohn(-Z)$ is an integral lattice for $(\oohm,\ang{y})$. For example,  $\oohn(-Z)(\usith)=w\cdot\cx[z,w]$, and the image of the bilinear map
$$w\cdot\cx[z,w]\times w\cdot\cx[z,w]\rf{z^n/w}\cx[z,w,1/z,1/w]$$
lies in $\V^0_{\oohn}(\usith)=\cx[z,w]$.  To find its dual lattice, note that
\bea{\oohn(-Z)'(\usith)&:=&\setst{f\in\cx[z,w,1/z,1/w]}{f\cdot \frac{z^n}{w}\cdot w\cdot\cx[z,w]\sst\cx[z,w]}\\
&=&\frac{1}{z^n}\cdot\cx[z,w]\\
&=& \oohn(X+nW)(\usith).
} 
We can similarly check on the other affine open subsets to conclude that $\oohn(-Z)'=\oohn(X+nW)$. Hence, $\Ll^0_{H^1_n}(\ang{y})$ is given by
\beq{\frac{\oo_{H_n}(X+nW)}{\oo_{H_n}(-Z)}\times\frac{\oo_{H_n}(X+nW)}{\oo_{H_n}(-Z)}\rf{\ang{y}}\frac{i_*\oo_{H_n-H^1_n}}{\oo_{H_n}}.\label{m1}}
Now we apply the map
$$\K^1_{H^1_n}:W(\cm^1_{H^1_n}(H_n))\ra W(\cm^1_{H^1_n-N}(H_n-N))\coprod W(\cm^1_{H^1_n-L}(H_n-L))$$
by restricting the domains to $H_n-N$ and $H_n-L$.
Restricting (\ref{m1}) to $H_n-L$ gives
\beq{\frac{\oo_{H_n-L}(nT_W)}{\oo_{H_n-L}}\times\frac{\oo_{H_n-L}(nT_W)}{\oo_{H_n-L}}\rf{\ang{y}}\frac{i\ds\oo_{H_n-H^1_n}}{\oo_{H_n-L}},\label{h0d0y1}}
while restricting it to $H_n-N$ gives
\beq{\frac{\oo_{H_n-N}(T_X)}{\oo_{H_n-N}(-T_Z)}\times\frac{\oo_{H_n-N}(T_X)}{\oo_{H_n-N}(-T_Z)}\rf{\ang{y}}\frac{i\ds\oo_{H_n-H^1_n}}{\oo_{H_n-N}}.\label{h0d0y2}}
On $(H_n-L)\cap\usio=\spec\cx[x,y,1/y]$, (\ref{h0d0y1}) is zero because
$$\left.\frac{\oo_{H_n-L}(nT_W)}{\oo_{H_n-L}}\right|_{(\hn-L)\cap\usio}=\frac{\oo_{(\hn-L)\cap\usio}}{\oo_{(\hn-L)\cap\usio}}=0,$$
while on $(H_n-L)\cap\usith=\spec\cx[z,w,1/w]$, it becomes
\beq{\frac{\frac{1}{z^n}\cdot\cx[z,w,1/w]}{\cx[z,w,1/w]}\times\frac{\frac{1}{z^n}\cdot\cx[z,w,1/w]}{\cx[z,w,1/w]} \rf{z^n/w}  \frac{\cx[z,w,1/z,1/w]}{\cx[z,w,1/w]}.\label{ann}}
If $n$ is even, then 
$$\frac{\frac{1}{z^{n/2}}\cdot\cx[z,w,1/w]}{z\cdot\cx[z,w,1/w]}\sst\frac{\frac{1}{z^n}\cdot\cx[z,w,1/w]}{z\cdot\cx[z,w,1/w]}$$
is a totally isotropic subspace of (\ref{ann}) of half the rank, i.e., a sublagrangian. Hence, the Witt class of (\ref{ann}) is zero, i.e., $\K^1_L\cc\Ll^0_{H^1_{\mb{\scriptsize even}}}(\ang{y})=0$. This implies that $d^0_{H_{\mb{\scriptsize even}}}(\ang{y})$ has no component in the subspace spanned by $\ang{1_y},\ang{y},\ang{1_w},\ang{w}$, and the corresponding rows in the third column are zero.

On the other hand, if $n$ is odd, then 
$$M:=\frac{\frac{1}{z^{(n-1)/2}}\cdot\cx[z,w,1/w]}{\cx[z,w,1/w]}\sst\frac{\frac{1}{z^n}\cdot\cx[z,w,1/w]}{\cx[z,w,1/w]}$$
is a totally isotropic subspace of (\ref{ann}), and
$$M\pp=\frac{\frac{1}{z^{(n+1)/2}}\cdot\cx[z,w,1/w]}{\cx[z,w,1/w]},$$
so $M\pp/M\simeq\cx[w,1/w]$, and \tcd
\beq{\xymatrix{\frac{\frac{1}{z^{(n+1)/2}}\cdot\cx[z,w,1/w]}{\frac{1}{z^{(n-1)/2}}\cdot\cx[z,w,1/w]}\times\frac{\frac{1}{z^{(n+1)/2}}\cdot\cx[z,w,1/w]}{\frac{1}{z^{(n-1)/2}}\cdot\cx[z,w,1/w]} \ar[r]^(.7){z^n/w} & \frac{\cx[z,w,1/z,1/w]}{\cx[z,w,1/w]}\\
\ar[u]_(.4){\rotatebox{90}{$\sim$}}^(.4){1/z^{(n+1)/2}\times1/z^{(n+1)/2}} \cx[w,1/w]\times\cx[w,1/w] \ar[r]^(.63){1/w} & \cx[w,1/w] \ar@{^(->}[u]_{1/z}
} \label{ann2}}
By Lemma \ref{final}(2) and Lemma \ref{oba2}, $M\mt M\pp/M$ does not change the Witt class, so we conclude that $\K^1_L\cc\Ll^0_{H^1_{\scriptsize odd}}(\ang{y})=\ang{1/w}=\ang{w}$. Hence, in the third column, the row corresponding to $\ang{w}$ is 1, while the rows corresponding to $\ang{1_y},\ang{y},\ang{1_w}$ are 0.

Now on $(\hn-N)\cap\usio=\spec\cx[x,y,1/x]$, (\ref{h0d0y2}) gives a \cd
\beq{\xymatrix{\frac{\frac{1}{y}\cdot\cx[x,y,1/x]}{\cx[x,y,1/x]}\times\frac{\frac{1}{y}\cdot\cx[x,y,1/x]}{\cx[x,y,1/x]} \ar[r]^(.63){y} & \frac{\cx[x,y,1/x,1/y]}{\cx[x,y,1/x]}\\
\dssimu^{1/y\times1/y} \cx[x,1/x]\times\cx[x,1/x] \ar[r]^(.63)1 & \cx[x,1/x] \ar@{^(->}[u]_{1/y}
}\label{head3}}
while on $(\hn-N)\cap\usith=\spec\cx[z,w,1/z]$, it gives a \cd
\beq{\xymatrix{\frac{\cx[z,w,1/z]}{w\cdot\cx[z,w,1/z]}\times\frac{\cx[z,w,1/z]}{w\cdot\cx[z,w,1/z]} \ar[r]^(.6){z^n/w} & \frac{\cx[z,w,1/z,1/w]}{\cx[z,w,1/z]}\\
\dssimu^{1\times1} \cx[z,1/z]\times\cx[z,1/z] \ar[r]^(.62){z^n} & \cx[z,1/z] \ar@{^(->}[u]_{1/w}
}\label{head4}}
Hence, we conclude that $\K^1_N\cc\Ll^0_{H^1_n}(\ang{y})=\ang{1_x}+\ang{z^n}$. 
This implies that on the third column, the rows corresponding to $\ang{1_x}$ and $\ang{z^n}$ are 1, while the rows corresponding to $\ang{x}$ and $\ang{z^{n+1}}$ are 0.

Finally, let us determine the entries in the fourth column. The form $(\oohm,\ang{xy})$ is given on the affine charts by
\bea{U_{\sigma_1} \quad :\quad \cx[x,y,1/x,1/y]\times\cx[x,y,1/x,1/y]&\rf{xy}&\cx[x,y,1/x,1/y]\\
U_{\sigma_2} \quad :\quad\cx[x,\by,1/x,1/\by]\times\cx[x,\by,1/x,1/\by]&\rf{x/\by}&\cx[x,\by,1/x,1/\by]\\
U_{\sigma_3} \quad :\quad\cx[z,w,1/z,1/w]\times\cx[z,w,,1/z,1/w]&\rf{z^{n-1}/w}&\cx[z,w,1/z,1/w]\\
U_{\sigma_4} \quad :\quad\cx[z,\bw,1/z,1/\bw]\times\cx[z,\bw,1/z,1/\bw]&\rf{z^{n-1}\bw}&\cx[z,\bw,1/z,1/\bw]
}
Using the same argument as above, we conclude that $\oo_{H_n}(-Z-W)$ is an integral lattice. For example, $\oohn(-Z-W)(\usith)=zw\cdot\cx[z,w]$, and the image of the bilinear form
$$zw\cdot\cx[z,w]\times zw\cdot\cx[z,w]\rf{z^{n-1}/w}\cx[z,w,1/z,1/w]$$
lies in $\cx[z,w]$.
To find its dual lattice, note that
\bea{\oohn(-Z-W)'(\usith)&:=&\setst{f\in\cx[z,w,1/z,1/w]}{f\cdot \frac{z^{n-1}}{w}\cdot zw\cdot\cx[z,w]\sst\cx[z,w]}\\
&=&\frac{1}{z^n}\cdot\cx[z,w]\\
&=& \oohn(X+Y+nW)(\usith).
} 
We can similarly check on the other affine open subsets to conclude that $\oohn(-Z-W)'=\oohn(X+Y+nW)$. Hence, $\Ll^0_{H^1_n}(\ang{xy})$ is given by
\beq{\frac{\oo_{H_n}(X+Y+nW)}{\oo_{H_n}(-Z-W)}\times\frac{\oo_{H_n}(X+Y+nW)}{\oo_{H_n}(-Z-W)}\rf{\ang{xy}}\frac{i\ds\oo_{H_n-H^1_n}}{\oo_{H_n}}.\label{head}}
Now we apply the map
$$\K^1_{H^1_n}:W(\cm^1_{H^1_n}(H_n))\ra W(\cm^1_{H^1_n-N}(H_n-N))\coprod W(\cm^1_{H^1_n-L}(H_n-L))$$
by restricting the domains to $H_n-N$ and $H_n-L$.

Restricting (\ref{head}) to $H_n-N$ gives
\beq{\frac{\oo_{H_n-N}(T_X)}{\oo_{H_n-N}(-T_Z)}\times\frac{\oo_{H_n-N}(T_X)}{\oo_{H_n-N}(-T_Z)}\rf{xy}\frac{i\ds\oo_{H_n-H^1_n}}{\oo_{H_n-N}},\label{h0d0xy1}}
while restricting to $H_n-L$ gives
\beq{\frac{\oo_{H_n-L}(T_Y+nT_W)}{\oo_{H_n-L}(-T_W)}\times\frac{\oo_{H_n-L}(T_Y+nT_W)}{\oo_{H_n-L}(-T_W)}\rf{xy}\frac{i\ds\oo_{H_n-H^1_n}}{\oo_{H_n-L}}.\label{h0d0xy2}}
On $(\hn-N)\cap\usio=\spec\cx[x,y,1/x]$, (\ref{h0d0xy1}) gives a \cd
\xym{\frac{\frac{1}{y}\cdot\cx[x,y,1/x]}{\cx[x,y,1/x]}\times\frac{\frac{1}{y}\cdot\cx[x,y,1/x]}{\cx[x,y,1/x]} \ar[r]^(.62){xy} & \frac{\cx[x,y,1/x,1/y]}{\cx[x,y,1/x]}\\
\dssimu^{1/y\times1/y} \cx[x,1/x]\times\cx[x,1/x] \ar[r]^(.62)x & \cx[x,1/x] \ar@{^(->}[u]_{1/y}
} while on $(\hn-N)\cap\usith=\spec\cx[z,w,1/z]$, it gives a \cd
\xym{\frac{\cx[z,w,1/z]}{w\cdot\cx[z,w,1/z]}\times\frac{\cx[z,w,1/z]}{w\cdot\cx[z,w,1/z]} \ar[r]^(.6){z^{n-1}/w} & \frac{\cx[z,w,1/z,1/w]}{\cx[z,w,1/z]}\\
\dssimu^{1\times1} \cx[z,1/z]\times\cx[z,1/z] \ar[r]^(.63){z^{n -1}} & \cx[z,1/z] \ar@{^(->}[u]_{1/w}
} Hence, we conclude that $\K^1_N\cc\Ll^0_{H^1_n}(\ang{xy})=\ang{x}+\ang{z^{n-1}}$. This implies that in the fourth column, the rows corresponding to $\ang{x}$ and $\ang{z^{n-1}}$ are 1, while the rows corresponding to $\ang{1_x}$ and $\ang{z^n}$ are 0.

On the other hand, on $(\hn-L)\cap\usio=\spec\cx[x,y,1/y]$, (\ref{h0d0xy2}) gives a \cd
\beq{\xymatrix{\frac{\frac{1}{x}\cdot\cx[x,y,1/y]}{\cx[x,y,1/y]}\times\frac{\frac{1}{x}\cdot\cx[x,y,1/y]}{\cx[x,y,1/y]} \ar[r]^(.62){xy} & \frac{\cx[x,y,1/x,1/y]}{\cx[x,y,1/y]}\\
\dssimu^{1/x\times1/x} \cx[y,1/y]\times\cx[y,1/y] \ar[r]^(.61)y & \cx[y,1/y] \ar@{^(->}[u]_{1/x}
}\label{head2}}
while on $(\hn-L)\cap\usith=\spec\cx[z,w,1/w]$, it becomes
\beq{\frac{\frac{1}{z^n}\cdot\cx[z,w,1/w]}{z\cdot\cx[z,w,1/w]}\times\frac{\frac{1}{z^n}\cdot\cx[z,w,1/w]}{z\cdot\cx[z,w,1/w]} \rf{z^{n-1}/w}  \frac{\cx[z,w,1/z,1/w]}{\cx[z,w,1/w]}.\label{tttt}}
If $n$ is odd, then
$$\frac{\frac{1}{z^{(n-1)/2}}\cdot\cx[z,w,1/w]}{z\cdot\cx[z,w,1/x]}\sst\frac{\frac{1}{z^n}\cdot\cx[z,w,1/w]}{z\cdot\cx[z,w,1/x]}$$
is a totally isotropic subspace of (\ref{tttt}) of half the rank, i.e., a sublagrangian. Hence, the Witt class of (\ref{tttt}) is zero, and together with (\ref{head2}), we conclude that
$$\K^1_{L}\cc\Ll^0_{H^1_{\mb{\scriptsize odd}}}(\ang{xy})=\ang{y}.$$
This implies that in the fourth column, the row corresponding to $\ang{y}$ is 1, while the rows corresponding to $\ang{1_y},\ang{1_w},\ang{w}$ are 0.

On the other hand, if $n$ is even, then 
$$M:=\frac{\frac{1}{z^{n/2-1}}\cdot\cx[z,w,1/w]}{z\cdot\cx[z,w,1/w]}\sst\frac{\frac{1}{z^n}\cdot\cx[z,w,1/w]}{z\cdot\cx[z,w,1/w]}$$
is a totally isotropic subspace of (\ref{tttt}), and
$$M\pp=\frac{\frac{1}{z^{n/2}}\cdot\cx[z,w,1/w]}{z\cdot\cx[z,w,1/w]},$$
so $M\pp/M\simeq\cx[w,1/w]$, and \tcd
\beq{\xymatrix{\frac{\frac{1}{z^{n/2}}\cdot\cx[z,w,1/w]}{\frac{1}{z^{n/2-1}}\cdot\cx[z,w,1/w]}\times\frac{\frac{1}{z^{n/2}}\cdot\cx[z,w,1/w]}{\frac{1}{z^{n/2-1}}\cdot\cx[z,w,1/w]} \ar[r]^(.67){z^{n-1}/w} & \frac{\cx[z,w,1/z,1/w]}{\cx[z,w,1/w]}\\
\ar[u]_(.4){\rotatebox{90}{$\sim$}}^(.4){1/z^{n/2}\times1/z^{n/2}} \cx[w,1/w]\times\cx[w,1/w] \ar[r]^(.63){1/w} & \cx[w,1/w] \ar@{^(->}[u]_{1/z}
}\label{t5}} 
As we noted earlier, $M\mt M\pp/M$ does not change the Witt class. Hence, together with (\ref{head2}), we conclude that 
$$\K^1_L\cc\Ll^1_{H^1_{\mb{\scriptsize even}}}(\ang{xy})=\ang{y}+\ang{1/w}=\ang{y}+\ang{w}.$$
This implies that on the fourth column, the rows corresponding to $\ang{y}$ and $\ang{w}$ are 1, while the rows corresponding to $\ang{1_y}$ and $\ang{1_w}$ are 0.
This completes the proof.
}

\prop{With the same choice of basis as in Proposition \ref{main1}, the map $d^1_{H_n}$ is represented by the matrix
\[
d^1_{H_{\mb{\scriptsize even}}}\quad=\quad
  \begin{blockarray}{ccccccccc}
    & \ang{1_x} & \ang{x} & \ang{1_y} & \ang{y} & \ang{1_z} & \ang{z} & \ang{1_w} & \ang{w} \\
    \begin{block}{c(cccccccc)}
   \ang{0_{xy}} & 0 & 1 & 0 & 1 & 0 & 0 & 0 & 0 \\
\ang{0_{ z\bw}} & 0 & 1 & 0 & 0 & 0 & 0 & 0 & 1  \\
 \ang{0_{x\by}} & 0 & 0 & 0 & 1 & 0 & 1 & 0 & 0 \\
\ang{0_{zw}}    & 0 & 0 & 0 & 0 & 0 & 1 & 0 & 1  \\
    \end{block}
      \end{blockarray}\quad,
\]
\[
d^1_{H_{\mb{\scriptsize odd}}}\quad=\quad
  \begin{blockarray}{ccccccccc}
    & \ang{1_x} & \ang{x} & \ang{1_y} & \ang{y} & \ang{1_z} & \ang{z} & \ang{1_w} & \ang{w} \\
    \begin{block}{c(cccccccc)}
   \ang{0_{xy}} & 0 & 1 & 0 & 1 & 0 & 0 & 0 & 0 \\
\ang{0_{ z\bw}} & 1 & 0 & 0 & 0 & 0 & 0 & 0 & 1  \\
 \ang{0_{x\by}} & 0 & 0 & 0 & 1 & 1 & 0 & 0 & 0 \\
\ang{0_{zw}}    & 0 & 0 & 0 & 0 & 0 & 1 & 0 & 1  \\
    \end{block}
      \end{blockarray}\quad.
\]
}

\pf{Recall that $d^1_{\hn}:=\Ll^1_{T_L}\coprod\Ll^1_{T_N}$, where
\bea{\Ll^1_{T_L}:W(\cm^1_{T_L}(\hn-N))&\ra& W(\cm^2_{\hn^2}(\hn)),\\
\Ll^1_{T_N}:W(\cm^1_{T_N}(\hn-L))&\ra& W(\cm^2_{\hn^2}(\hn)).
}
$W(\cm^1_{T_L}(H_n-N))$ is generated by $\ang{1_x}, \ang{x}, \ang{1_z}, \ang{z}$, and $W(\cm^1_{T_N}(H_n-L))$ is generated by $\ang{1_y}, \ang{y}, \ang{1_w}, \ang{w}$ (Corollary \ref{ttt}).
\Tcd
\xym{W(\cm^1_{T_L}(H_n-N)) \ar[r]^(.55){\Ll^1_{T_L}} & W(\cm^2_{H^2}(H_n)) \\
W(\cm^1_{T_X}(\usio-Y)) \ar[u] \ar[r]^(.55){\Ll^1_{\sigma_1, y}}  & \ar[u] W(\cm^2_{0_{xy}}(\usio))
}  where the vertical maps are induced by inclusion. Hence, we can compute $\Ll^1_{T_L}$ in terms of the affine lattice maps $\Ll^1_{\sigma_1,y}, \Ll^1_{\sigma_2,\by}, \Ll^1_{\sigma_3,w}, \Ll^1_{\sigma_4,\bw}$, and $\Ll^1_{T_N}$ in terms of $\Ll^1_{\sigma_1,x}, \Ll^1_{\sigma_2,x}, \Ll^1_{\sigma_3,z}, \Ll^1_{\sigma_4,z}$. The integral lattice and dual lattice for $\Ll^1_{\sigma_1,y}$ are computed using the value group
\beq{\V^1_{\oo_{\usio}, y}=\frac{\cx[x,y,1/y]}{\cx[x,y]}\simeq\frac{\cx[x,y,1/y]+\cx[x,y,1/x]}{\cx[x,y,1/x]}\sst\frac{\cx[x,y,1/x,1/y]}{\cx[x,y,1/x]}.\label{sun3}}
We saw in the computation of $\Ll^0_{H_n^1}(\ang{y})$  that  on $(\hn-N)\cap\usio=\spec\cx[x,y,1/x]$, $\ang{1_x}$ is given by (\ref{head3})
\beq{\frac{\frac{1}{y}\cdot\cx[x,y,1/x]}{\cx[x,y,1/x]}\times\frac{\frac{1}{y}\cdot\cx[x,y,1/x]}{\cx[x,y,1/x]} \rf{y} \frac{\cx[x,y,1/x,1/y]}{\cx[x,y,1/x]}.\label{abou}}
Let $M:=\frac{\frac{1}{y}\cdot\cx[x,y]+\cx[x,y,1/x]}{\cx[x,y,1/x]}\sst\frac{\frac{1}{y}\cdot\cx[x,y,1/x]}{\cx[x,y,1/x]}$. Since the image of the bilinear form
$$\frac{\frac{1}{y}\cdot\cx[x,y]+\cx[x,y,1/x]}{\cx[x,y,1/x]}\times\frac{\frac{1}{y}\cdot\cx[x,y]+\cx[x,y,1/x]}{\cx[x,y,1/x]} \rf{y} \frac{\cx[x,y,1/x,1/y]}{\cx[x,y,1/x]}$$
lies in $\V^1_{\oo_{\usio},y}$ (\ref{sun3}), $M$ is an integral lattice for (\ref{abou}). Moreover,
$$M'=\setst{f\in\frac{\frac{1}{y}\cdot\cx[x,y,1/x]}{\cx[x,y,1/x]}}{f\cdot \frac{1}{y}\cdot y\in\frac{\frac{1}{y}\cdot\cx[x,y]+\cx[x,y,1/x]}{\cx[x,y,1/x]}}=M,$$
so $M$ is self-dual. Hence, $\Ll^1_{\sigma_1, y}(\ang{1_x})=0$. This implies that $\Ll^1_{T_L}(\ang{1_x})$ has no $\ang{0_{xy}}$ component, so the corresponding row in the first column is 0.

On $(\hn-N)\cap\usif=\spec[z,\bw,1/z]$, $\ang{1_x}$ is represented by (see \ref{abou})
\beq{\frac{\frac{1}{\bw}\cdot\cx[ z,\bw,1/ z]}{\cx[ z,\bw,1/ z]}\times\frac{\frac{1}{\bw}\cdot\cx[ z,\bw,1/ z]}{\cx[ z,\bw,1/ z]} \rf{z^n\bw} \frac{\cx[ z,\bw,1/z,1/\bw]}{\cx[ z,\bw,1/ z]}.\label{nigh}}
Using the the value group
\beq{\V^1_{\oo_{\usif}, \bw}=\frac{\cx[z,\bw,1/\bw]}{\cx[z,\bw]}\simeq\frac{\cx[z,\bw,1/\bw]+\cx[z,\bw,1/z]}{\cx[z,\bw,1/z]}\sst\frac{\cx[z,\bw,1/z,1/\bw]}{\cx[z,\bw,1/z]},\label{sun}}
one finds that $M:=\frac{\frac{1}{\bw}\cdot\cx[ z,\bw]+\cx[ z,\bw,1/ z]}{\cx[ z,\bw,1/ z]}$ is an integral lattice, and $M'=\frac{\frac{1}{z^n\bw}\cdot\cx[ z,\bw]+\cx[ z,\bw,1/ z]}{\cx[ z,\bw,1/ z]}$ is its dual lattice. 
Since the rank \hmm\ induces
$$W(\cm^2_{0_{xy}}(\usio))=W(\cx)=\Z/2,$$
$\dm_\cx(M'/M)=n$ implies that  $\Ll^1_{\sigma_4, \bw}(\ang{1_x})$ is zero \ifof $n$ is even. This implies that $\Ll^1_{T_L}(\ang{1_x})$ has no $\ang{0_{z\bw}}$ component \ifof $n$ is even, so the corresponding row in the first column is zero \ifof $n$ is even.

Since $\ang{1_x}$ is supported on $T_X$ (Corollary \ref{ttt}), the rest of the entries in the first column are zero.

The entries in the second column are obtained in the same way. 
On $(\hn-N)\cap\usio=\spec\cx[x,y,1/x]$, $\ang{x}$ is represented by (see \ref{abou})
\beq{\frac{\frac{1}{y}\cdot\cx[x,y,1/x]}{\cx[x,y,1/x]}\times\frac{\frac{1}{y}\cdot\cx[x,y,1/x]}{\cx[x,y,1/x]} \rf{xy} \frac{\cx[x,y,1/x,1/y]}{\cx[x,y,1/x]}.\label{abou2}}
One finds that $M=\frac{\frac{1}{y}\cdot\cx[x,y]+\cx[x,y,1/x]}{\cx[x,y,1/x]}$ is an integral lattice, but and that its dual lattice is $M'=\frac{\frac{1}{xy}\cdot\cx[x,y]+\cx[x,y,1/x]}{\cx[x,y,1/x]}$. Then $\dm_\cx M'/M=1$, so that $\Ll^1_{\sigma_1, y}(\ang{x})$ is non-zero. This implies that $\Ll^1_{T_L}(\ang{x})$ has a non-zero $\ang{0_{xy}}$ component, so the corresponding row in the second column is 1.

On $(\hn-N)\cap\usif=\spec[z,\bw,1/z]$, $\ang{x}$ is represented by (see \ref{abou2})
\beq{\frac{\frac{1}{\bw}\cdot\cx[ z,\bw,1/ z]}{\cx[ z,\bw,1/ z]}\times\frac{\frac{1}{\bw}\cdot\cx[ z,\bw,1/ z]}{\cx[ z,\bw,1/ z]} \rf{z^{n-1}\bw} \frac{\cx[ z,\bw,1/z,1/\bw]}{\cx[ z,\bw,1/ z]}.\label{nigh2}}
Using the value group $\V^1_{\oo_{\usif}, \bw}$ (\ref{sun}), one finds that $M:=\frac{\frac{z}{\bw}\cdot\cx[ z,\bw]+\cx[ z,\bw,1/ z]}{\cx[ z,\bw,1/ z]}$ is an integral lattice, and $M'=\frac{\frac{1}{z^n\bw}\cdot\cx[ z,\bw]+\cx[ z,\bw,1/ z]}{\cx[ z,\bw,1/ z]}$ is its dual lattice. Then $\dm_\cx M'/M=n+1$, so $\Ll^1_{\sigma_4,\bw}(\ang{x})$ is zero \ifof $n$ is odd. This implies that $\Ll^1_{T_L}(\ang{x})$ has no $\ang{0_{z\bw}}$ component \ifof $n$ is odd, so the corresponding row in the second column is 0 \ifof $n$ is odd.

Since $\ang{x}$ is supported on $T_X$, the rest of the entries in the second column are zero.

We move on to the third column.

We saw in the computation of $\Ll^0_{H^1_n}(\ang{x})$ that on $(\hn-L)\cap\usio=\spec\cx[x,y,1/y]$, $\ang{1_y}$ is given by (\ref{sun2})
\beq{\frac{\frac{1}{x}\cdot\cx[x,y,1/y]}{\cx[x,y,1/y]}\times\frac{\frac{1}{x}\cdot\cx[x,y,1/y]}{\cx[x,y,1/y]} \rf{x} \frac{\cx[x,y,1/x,1/y]}{\cx[x,y,1/y]}.\label{sun4}}
Using the value group 
\beq{\V^1_{\oo_{\usio}, x}=\frac{\cx[x,y,1/x]}{\cx[x,y]}\simeq\frac{\cx[x,y,1/x]+\cx[x,y,1/y]}{\cx[x,y,1/y]}\sst\frac{\cx[x,y,1/x,1/y]}{\cx[x,y,1/y]},\label{valo}}
one finds that $\frac{\frac{1}{x}\cdot\cx[x,y]+\cx[x,y,1/y]}{\cx[x,y,1/y]}$ is a self-dual integral lattice, so $\Ll^1_{\sigma_1, x}(\ang{1_y})=0$. This implies that $\Ll^1_{T_N}(\ang{1_y})$ has no $\ang{0_{xy}}$ component, so the corresponding row in the third column is 0.

On $(\hn-L)\cap\usit=\spec\cx[x,\by,1/\by]$, $\ang{1_y}$ is given by (see (\ref{sun4}))
\beq{\frac{\frac{1}{x}\cdot\cx[x,\by,1/\by]}{\cx[x,y,1/y]}\times\frac{\frac{1}{x}\cdot\cx[x,\by,1/\by]}{\cx[x,\by,1/\by]} \rf{x} \frac{\cx[x,\by,1/x,1/\by]}{\cx[x,\by,1/\by]}.\label{sun6}}
Using the value group
\beq{\V^1_{\oo_{\usit}, x}=\frac{\cx[x,\by,1/x]}{\cx[x,\by]}\simeq\frac{\cx[x,\by,1/x]+\cx[x,\by,1/\by]}{\cx[x,\by,1/\by]}\sst\frac{\cx[x,\by,1/x,1/\by]}{\cx[x,\by,1/\by]},
\label{sun7}}
one finds that $\frac{\frac{1}{x}\cdot\cx[x,\by]+\cx[x,\by,1/\by]}{\cx[x,\by,1/\by]}$ is a self-dual integral lattice, so $\Ll^1_{\sigma_2, x}(\ang{1_y})=0$. This implies that $\Ll^1_{T_N}(\ang{1_y})$ has no $\ang{0_{x\by}}$ component, so the corresponding row in the third column is 0.

Since $\ang{1_y}$ is supported on $T_Y$, the rest of the entries in the third column are zero.

We move on to the fourth column.

On $(\hn-L)\cap\usio$, $\ang{y}$ is represented by (see (\ref{sun4}))
\beq{\frac{\frac{1}{x}\cdot\cx[x,y,1/y]}{\cx[x,y,1/y]}\times\frac{\frac{1}{x}\cdot\cx[x,y,1/y]}{\cx[x,y,1/y]} \rf{xy} \frac{\cx[x,y,1/x,1/y]}{\cx[x,y,1/y]}.
\label{beau}} 
Using the value group $\V^1_{\oo_{\usio},x}$ (\ref{valo}), one finds that $M:=\frac{\frac{1}{x}\cdot\cx[x,y]+\cx[x,y,1/y]}{\cx[x,y,1/y]}$ is an integral lattice, and that $M':=\frac{\frac{1}{xy}\cdot\cx[x,y]+\cx[x,y,1/y]}{\cx[x,y,1/y]}$ is its dual lattice. Since $\dm_\cx M'/M=1$, we conclude that $\Ll^1_{\sigma_2, x}(\ang{y})\neq0$. This implies that $\Ll^1_{T_N}(\ang{y})$ has a nonzero $\ang{0_{xy}}$ component, so the corresponding row in the fourth column is 1.

On $(\hn-L)\cap\usit$, $\ang{y}$ is represented by (see \ref{sun6})
\beq{\frac{\frac{1}{x}\cdot\cx[x,\by,1/\by]}{\cx[x,y,1/y]}\times\frac{\frac{1}{x}\cdot\cx[x,\by,1/\by]}{\cx[x,\by,1/\by]} \rf{x/\by} \frac{\cx[x,\by,1/x,1/\by]}{\cx[x,\by,1/\by]}\label{beau2}
} 
Using the value group $\V^1_{\oo_{\usit},x}$ (\ref{sun7}), one finds that $M:=\frac{\frac{\by}{x}\cdot\cx[x,\by]+\cx[x,\by,1/\by]}{\cx[x,\by,1/\by]}$ is an integral lattice, and that $M':=\frac{\frac{1}{x}\cdot\cx[x,\by]+\cx[x,\by,1/\by]}{\cx[x,\by,1/\by]}$ is its dual lattice. Then $\dm_\cx M'/M=1$, so $\Ll^1_{\sigma_2, x}(\ang{y})\neq0$. This implies that $\Ll^1_{T_N}(\ang{y})$ has a nonzero $\ang{0_{x\by}}$ component, so the corresponding row in the fourth column is 1.

Since $\ang{y}$ is supported on $T_Y$, the rest of the entries in the fourth column are zero.

We move on to the fifth column. 

We saw in the computation of $\Ll^0_{\hn^1}(\ang{y})$ that on $(\hn-N)\cap\usith=\cx[z,w,1/z]$, $\ang{1_z}$ is represented by (\ref{head4})
\beq{\frac{\cx[z,w,1/z]}{w\cdot\cx[z,w,1/z]}\times\frac{\cx[z,w,1/z]}{w\cdot\cx[z,w,1/z]} \rf{1/w} \frac{\cx[z,w,1/z,1/w]}{\cx[z,w,1/z]}.\label{head5}}
Using the value group 
\beq{\V^1_{\oo_{\usith}, w}=\frac{\cx[z,w,1/w]}{\cx[z,w]}\simeq\frac{\cx[z,w,1/w]+\cx[z,w,1/z]}{\cx[z,w,1/z]}\sst\frac{\cx[z,w,1/z,1/w]}{\cx[z,w,1/z]},\label{sun5}}
one finds that $\frac{\cx[z,w]+w\cdot\cx[z,w,1/x]}{w\cdot\cx[z,w,1/x]}$ is a self-dual integral lattice, so $\Ll^1_{\sigma_3, w}(\ang{1_z})=0$. This implies that $\Ll^1_{T_L}(\ang{1_z})$ has no $\ang{0_{zw}}$ component, so the corresponding row in the fifth column is 0.

On $(\hn-N)\cap\usit=\cx[x,\by,1/x]$, $\ang{1_z}$ is represented  by (see \ref{head5})
\beq{\frac{\cx[x,\by,1/x]}{\by\cdot\cx[x,\by,1/x]}\times\frac{\cx[x,\by,1/x]}{\by\cdot\cx[x,\by,1/x]} \rf{x^n/\by} \frac{\cx[x,\by,1/x,1/\by]}{\cx[x,\by,1/x]},\label{chri}}
and using the value group $\V^1_{\oo_{\usit}, x}$ (\ref{sun7}), one finds that $M:=\frac{\cx[x,\by]+\by\cdot\cx[x,\by,1/x]}{\by\cdot\cx[x,\by,1/x]}$ is an  integral lattice, and that $M':=\frac{\frac{1}{x^n}\cdot\cx[x,\by]+\by\cdot\cx[x,\by,1/x]}{\by\cdot\cx[x,\by,1/x]}$ is its dual lattice. Hence, $\dm_\cx M'/M=n$, so $\Ll^1_{\sigma_2, x}(\ang{1_z})=0$ \ifof $n$ is even. This implies that $\Ll^1_{T_L}(\ang{1_z})$ has no $\ang{0_{x\by}}$ component \ifof $n$ is even, so the corresponding row in the fifth column is 0 \ifof $n$ is even.

Since $\ang{1_z}$ is supported on $T_Z$, the rest of the entries in the fifth column are zero.

We move on to the sixth column.

On $(\hn-N)\cap\usith=\cx[z,w,1/z]$, $\ang{z}$ is represented by (see \ref{head5})
\beq{\frac{\cx[z,w,1/z]}{w\cdot\cx[z,w,1/z]}\times\frac{\cx[z,w,1/z]}{w\cdot\cx[z,w,1/z]} \rf{z/w} \frac{\cx[z,w,1/z,1/w]}{\cx[z,w,1/z]}.
\label{head8}}
Using the value group $\V^1_{\oo_{\usith},w}$ (\ref{sun5}), one finds that $M:=\frac{\cx[z,w]+w\cdot\cx[z,w,1/z]}{w\cdot\cx[z,w,1/z]}$ is an integral lattice, and that $M':=\frac{\frac{1}{z}\cdot\cx[z,w]+w\cdot\cx[z,w,1/z]}{w\cdot\cx[z,w,1/z]}$ is its dual lattice. Hence, $\dm_\cx M'/M=1$, so $\Ll^1_{\sigma_3, w}(\ang{z})\neq0$. This implies that $\Ll^1_{T_L}(\ang{z})$ has a non-zero $\ang{0_{zw}}$ component, so the corresponding row in the sixth column is 1.

On $(\hn-N)\cap\usit=\cx[x,\by,1/x]$, $\ang{z}$ is represented by (see \ref{chri})
\beq{\frac{\cx[x,\by,1/x]}{\by\cdot\cx[x,\by,1/x]}\times\frac{\cx[x,\by,1/x]}{\by\cdot\cx[x,\by,1/x]} \rf{x^{n-1}/\by} \frac{\cx[x,\by,1/x,1/\by]}{\cx[x,\by,1/x]}.\label{head9}}

Using the value group 
\beq{\V^1_{\oo_{\usit},\by}=\frac{\cx[x,\by,1/\by]}{\cx[x,\by]}\simeq\frac{\cx[x,\by,1/\by]+\cx[x,\by,1/x]}{\cx[x,\by,1/x]}\sst\frac{\cx[x,\by,1/x,1/\by]}{\cx[x,\by,1/x]},
\label{val2by}}
one finds that $M=\frac{\cx[x,\by]+\by\cdot\cx[x,\by,1/x]}{\by\cdot\cx[x,\by,1/x]}$ is an integral lattice, and that $M'=\frac{\frac{1}{x^{n-1}}\cdot\cx[x,\by]+\by\cdot\cx[x,\by,1/x]}{\by\cdot\cx[x,\by,1/x]}$ is its dual lattice. Then $\dm_\cx M'/M=n-1$, so $\Ll^1_{\sigma_2, \by}(\ang{z})\neq0$ \ifof $n$ is even. This implies that $\Ll^1_{T_L}(\ang{z})$ has a non-zero $\ang{0_{x\by}}$ component \ifof $n$ is even, so the corresponding row in the sixth column is 1 \ifof $n$ is even.

Since $\ang{z}$ is supported on $T_Z$, the rest of the entries in the sixth column are zero.

For the seventh and eight columns, we consider different parities of $n$ separately at the outset:

(1) When $n$ is even : We saw in the computation of $\Ll^0_{\hn^1}(\ang{xy})$ that on $(\hn-L)\cap\usith=\spec\cx[z,w,1/w]$, $\ang{1_w}$ is given by (see (\ref{tttt}) and (\ref{t5}))
\beq{\frac{\frac{1}{z^n}\cdot\cx[z,w,1/w]}{z\cdot\cx[z,w,1/w]}\times\frac{\frac{1}{z^n}\cdot\cx[z,w,1/w]}{z\cdot\cx[z,w,1/w]} \rf{z^{n-1}/w^2} \frac{\cx[z,w,1/z,1/w]}{\cx[z,w,1/w]}.\label{toot}}
Using the value group
\beq{\V^1_{\oo_{\usith,z}}=\frac{\cx[z,w,1/z]}{\cx[z,w]}\simeq\frac{\cx[z,y,1/z]+\cx[z,y,1/w]}{\cx[z,w,1/w]}\sst\frac{\cx[z,w,1/z,1/w]}{\cx[z,w,1/w]},
\label{val3z}}
one finds that $\frac{\frac{w}{z^{n/2}}\cdot\cx[z,w]+z\cdot\cx[z,w,1/w]}{z\cdot\cx[z,w,1/w]}$  is a self-dual integral lattice, so $\Ll^1_{\sigma_3, z}(\ang{w})=0$. This implies that $\Ll^1_{T_N}(\ang{w})$ does not have   $\ang{0_{zw}}$ component, so the corresponding row  in the seventh column is 0.

On $(\hn-L)\cap\usif=\spec\cx[z,\bw,1/\bw]$, $\ang{1_w}$ is given by (see \ref{toot})
\beq{\frac{\frac{1}{z^n}\cdot\cx[z,\bw,1/\bw]}{z\cdot\cx[z,\bw,1/\bw]}\times\frac{\frac{1}{z^n}\cdot\cx[z,\bw,1/\bw]}{z\cdot\cx[z,\bw,1/\bw]} \rf{z^{n-1}\bw^2} \frac{\cx[z,\bw,1/z,1/\bw]}{\cx[z,\bw,1/\bw]}.\label{toot2}}
Using the value group 
\beq{\V^1_{\oo_{\usif}, z}=\frac{\cx[z,\bw,1/z]}{\cx[z,\bw]}\simeq\frac{\cx[z,\bw,1/z]+\cx[z,\bw,1/\bw]}{\cx[z,\bw,1/\bw]}\sst\frac{\cx[z,\bw,1/z,1/\bw]}{\cx[z,\bw,1/\bw]},\label{val4z}}
one finds that $\frac{\frac{1}{z^{n/2}\bw}\cdot\cx[z,\bw]+z\cdot\cx[z,\bw,1/\bw]}{z\cdot\cx[z,\bw,1/\bw]}$ is a self-dual integral lattice, so $\Ll^1_{\sigma_4,z}(\ang{1_w})=0$. This implies that $\Ll^1_{T_N}(\ang{1_w})$ has no $\ang{0_{z\bw}}$ component, so the corresponding row in the seventh column is 0.

Since $\ang{1_w}$ is supported on $T_W$, the rest of the entries in the seventh column are zero.

We move on to the eighth column.

On $(\hn-L)\cap\usith=\cx[z,w,1/w]$, $\ang{w}$ is represented by (see \ref{toot})
\beq{\frac{\frac{1}{z^n}\cdot\cx[z,w,1/w]}{z\cdot\cx[z,w,1/w]}\times\frac{\frac{1}{z^n}\cdot\cx[z,w,1/w]}{z\cdot\cx[z,w,1/w]} \rf{z^{n-1}/w} \frac{\cx[z,w,1/z,1/w]}{\cx[z,w,1/w]}.\label{toot3}}
Using the value group $\V^1_{\oo_{\usith,z}}$ (\ref{val3z}), one finds that $M:=\frac{\frac{w}{z^{n/2}}\cdot\cx[z,w]+z\cdot\cx[z,w,1/w]}{z\cdot\cx[z,w,1/w]}$ is an integral lattice, and that $M':=\frac{\frac{1}{z^{n/2}}\cdot\cx[z,w]+z\cdot\cx[z,w,1/w]}{z\cdot\cx[z,w,1/w]}$ is its dual lattice. Then $\dm_\cx M'/M=1$, so $\Ll^1_{\sigma_4, z}(\ang{w})\neq0$. This implies that $\Ll^1_{T_N}(\ang{w})$ has a non-zero $\ang{0_{zw}}$ component, so the corresponding row in the eighth column is 1.

On $(\hn-L)\cap\usif=\spec\cx[z,\bw,1/\bw]$, $\ang{w}$ is represented by (see \ref{toot2})
\beq{\frac{\frac{1}{z^n}\cdot\cx[z,\bw,1/\bw]}{z\cdot\cx[z,\bw,1/\bw]}\times\frac{\frac{1}{z^n}\cdot\cx[z,\bw,1/\bw]}{z\cdot\cx[z,\bw,1/\bw]} \rf{z^{n-1}\bw} \frac{\cx[z,\bw,1/z,1/\bw]}{\cx[z,\bw,1/\bw]}.\label{toot4}}
Using the value group $\V^1_{\oo_{\usif},z}$ (\ref{val4z}), one finds that $M:=\frac{\frac{1}{z^{n/2}}\cdot\cx[z,\bw]+z\cdot\cx[z,\bw,1/\bw]}{z\cdot\cx[z,\bw,1/\bw]}$ is an integral lattice, and that $M':=\frac{\frac{1}{z^{n/2}\bw}\cdot\cx[z,\bw]+z\cdot\cx[z,\bw,1/\bw]}{z\cdot\cx[z,\bw,1/\bw]}$ is its dual lattice. Hence, $\dm_\cx M'/M=1$, so $\Ll^1_{\sigma_4, z}(\ang{w})\neq0$. This implies that $\Ll^1_{T_N}(\ang{w})$ has a non-zero $\ang{0_{z\bw}}$ component, so the corresponding row in the eigth column is 1.

Since $\ang{w}$ is supported on $T_W$, the rest of the entries in the eighth column are zero.

(2) When $n$ is odd : We saw in the computation of $\Ll^0_{\hn^1}(\ang{y})$ that on $(\hn-L)\cap\usith=\cx[z,w,1/w]$, $\ang{1_w}$ is represented by (see (\ref{ann}) and (\ref{ann2}))
\beq{\frac{\frac{1}{z^n}\cdot\cx[z,w,1/w]}{\cx[z,w,1/w]}\times\frac{\frac{1}{z^n}\cdot\cx[z,w,1/w]}{\cx[z,w,1/w]} \rf{z^n/w^2} \frac{\cx[z,w,1/z,1/w]}{\cx[z,w,1/w]}.\label{hungry2}}
Using the value group $\V^1_{\oo_{\usith},z}$ (\ref{val3z}), one finds that $\frac{\frac{w}{z^{(n+1)/2}}\cdot\cx[z,w]+\cx[z,w,1/w]}{\cx[z,w,1/w]}$ is a self-dual integral lattice, so $\Ll^1_{\sigma_3, z}(\ang{w})=0$. This implies that $\Ll^1_{T_N}(\ang{w})$ has  no $\ang{0_{zw}}$ component, so the corresponding row in the seventh column is 0.

On $(\hn-L)\cap\usif=\spec\cx[z,\bw,1/\bw]$, $\ang{1_w}$ is given by (see (\ref{hungry2}))
\beq{\frac{\frac{1}{z^n}\cdot\cx[z,\bw,1/\bw]}{\cx[z,\bw,1/\bw]}\times\frac{\frac{1}{z^n}\cdot\cx[z,\bw,1/\bw]}{\cx[z,\bw,1/\bw]} \rf{z^n\bw^2} \frac{\cx[z,\bw,1/z,1/\bw]}{\cx[z,\bw,1/\bw]}.\label{hungry3}}
Using the value group $\V^1_{\oo_{\usif},z}$ (\ref{val4z}), one finds that $\frac{\frac{1}{z^{(n+1)/2}\bw}\cdot\cx[z,\bw]+\cx[z,\bw,1/\bw]}{\cx[z,\bw,1/\bw]}$ is a self-dual integral lattice, so $\Ll^1_{\sigma_4, z}(\ang{1_w})=0$. This implies that $\Ll^1_{T_N}(\ang{1_w})$ has  no $\ang{0_{z\bw}}$ component, so the corresponding row in the seventh column is 0.

On $(\hn-L)\cap\usith=\spec\cx[z,w,1/w]$, $\ang{w}$ is given by (see  (\ref{hungry2}))
\beq{\frac{\frac{1}{z^n}\cdot\cx[z,w,1/w]}{\cx[z,w,1/w]}\times\frac{\frac{1}{z^n}\cdot\cx[z,w,1/w]}{\cx[z,w,1/w]} \rf{z^n/w} \frac{\cx[z,w,1/z,1/w]}{\cx[z,w,1/w]}.\label{hungry}}
Using the value group $\V^1_{\oo_{\usith},z}$ (\ref{val3z}), one finds that $M:=\frac{\frac{w}{z^{(n+1)/2}}\cdot\cx[z,w]+\cx[z,w,1/w]}{\cx[z,w,1/w]}$ is an  integral lattice, and that $M':=\frac{\frac{1}{z^{(n+1)/2}}\cdot\cx[z,w]+\cx[z,w,1/w]}{\cx[z,w,1/w]}$ is its dual lattice. Hence, $\dm_\cx M'/M=1$, so $\Ll^1_{\sigma_4, z}(\ang{w})\neq0$. This implies that $\Ll^1_{T_N}(\ang{w})$ has a non-zero $\ang{0_{zw}}$ component, so the corresponding row in the eighth column is 1.

On $(\hn-L)\cap\usif=\spec\cx[z,\bw,1/\bw]$, $\ang{w}$ is given by (see (\ref{hungry3}))
\beq{\frac{\frac{1}{z^n}\cdot\cx[z,\bw,1/\bw]}{\cx[z,\bw,1/\bw]}\times\frac{\frac{1}{z^n}\cdot\cx[z,\bw,1/\bw]}{\cx[z,\bw,1/\bw]} \rf{z^n\bw} \frac{\cx[z,\bw,1/z,1/\bw]}{\cx[z,\bw,1/\bw]}.\label{hung2}}
Using the value group $\V^1_{\oo_{\usif},z}$ (\ref{val4z}), one finds that $M:=\frac{\frac{1}{z^{(n+1)/2}}\cdot\cx[z,\bw]+\cx[z,\bw,1/\bw]}{\cx[z,\bw,1/\bw]}$ is an integral lattice, and that $M':=\frac{\frac{1}{z^{(n+1)/2}\bw}\cdot\cx[z,\bw]+\cx[z,\bw,1/\bw]}{\cx[z,\bw,1/\bw]}$ is its dual lattice. Then $\dm_\cx M'/M=1$, so $\Ll^1_{\sigma_4, z}(\ang{w})\neq0$. This implies that $\Ll^1_{T_N}(\ang{w})$ has a non-zero $\ang{0_{z\bw}}$ component, so the corresponding row in the eighth column is 1.

Since $\ang{w}$ is supported on $T_W$, the rest of the entries in the eighth column are zero.
}

\rmk{The matrix representation for $d^1_{H^1_n}$ can also be deduced from Schmid's result \cite{schmid} that the canonical map 
$$W(\PP^1;\oo_{\PP^1}(n))\ra \coprod_{x\in \PP^{1(1)}}W(x)$$
is given by the second residue \hmm\ (Theorem \ref{springer}) at $x\neq\infty$, and by the first (resp. second) residue \hmm\ at $x=\infty$ if $n$ is odd (resp. even).
}
\chapter{Cohomologies}
\label{coho}

Now that we have seen the quasi-\iso\ between the toric complex (\ref{mon3}) and the Gersten-Witt complex  (\ref{mon2}) of $\hn$, we compute cohomologies using the former. We will see that they are cohomologies of the Witt sheaf $U\mt W(U)$ on $\hn$.

We first verify that $d^1_{H_n}\cc d^0_{H_n}=0$. 
Next, to compute the cohomologies, we put the matrices into Smith normal forms\footnote{They are computed using a Mathematica package \itk{IntegerSmithNormalForm} from  \url{http://library.wolfram.com/infocenter/MathSource/682}.} :
$$d^0_{H_{\mb{\scriptsize even}}}\simeq
  \begin{blockarray}{ccccc}
    & \ang{1} & \ang{x} & \ang{y} & \ang{xy} \\
    \begin{block}{c(cccc)}
    \ang{1_x} & 1 & 0 & 0 & 0\\
    \ang{x} & 0 & 1 & 0 & 0 \\
    \ang{1_y} & 0 & 0 & 1 & 0 \\
    \ang{y} & 0 & 0 & 0 & 0 \\
    \ang{1_z} & 0 & 0 & 0 & 0 \\
    \ang{z} & 0 & 0 & 0 & 0 \\
    \ang{1_w} & 0 & 0 & 0 & 0 \\
    \ang{w} & 0 & 0 & 0 & 0 \\
    \end{block}
      \end{blockarray} \quad,\quad
  d^0_{H_{\mb{\scriptsize odd}}}\simeq
  \begin{blockarray}{ccccc}
    & \ang{1} & \ang{x} & \ang{y} & \ang{xy} \\
    \begin{block}{c(cccc)}
    \ang{1_x} & 1 & 0 & 0 & 0\\
    \ang{x} & 0 & 1 & 0 & 0 \\
    \ang{1_y} & 0 & 0 & 1 & 0 \\
    \ang{y} & 0 & 0 & 0 & 0 \\
    \ang{1_z} & 0 & 0 & 0 & 0 \\
    \ang{z} & 0 & 0 & 0 & 0 \\
    \ang{1_w} & 0 & 0 & 0 & 0 \\
    \ang{w} & 0 & 0 & 0 & 0 \\
    \end{block}
      \end{blockarray}\quad.
$$
\[
d^1_{H_{\mb{\scriptsize even}}}\quad\simeq\quad
  \begin{blockarray}{ccccccccc}
    & \ang{1_x} & \ang{x} & \ang{1_y} & \ang{y} & \ang{1_z} & \ang{z} & \ang{1_w} & \ang{w} \\
    \begin{block}{c(cccccccc)}
   \ang{0_{xy}} & 1 & 0 & 0 & 0 & 0 & 0 & 0 & 0 \\
\ang{0_{ z\bw}} & 0 & 1 & 0 & 0 & 0 & 0 & 0 & 0  \\
 \ang{0_{x\by}} & 0 & 0 & 1 & 0 & 0 & 0 & 0 & 0 \\
\ang{0_{zw}}    & 0 & 0 & 0 & 0 & 0 & 0 & 0 & 0  \\
    \end{block}
      \end{blockarray}\quad,
\]
\[
d^1_{H_{\mb{\scriptsize odd}}}\quad\simeq\quad
  \begin{blockarray}{ccccccccc}
    & \ang{1_x} & \ang{x} & \ang{1_y} & \ang{y} & \ang{1_z} & \ang{z} & \ang{1_w} & \ang{w} \\
    \begin{block}{c(cccccccc)}
   \ang{0_{xy}} & 1 & 0 & 0 & 0 & 0 & 0 & 0 & 0 \\
\ang{0_{ z\bw}} & 0 & 1 & 0 & 0 & 0 & 0 & 0 & 0  \\
 \ang{0_{x\by}} & 0 & 0 & 1 & 0 & 0 & 0 & 0 & 0 \\
\ang{0_{zw}}    & 0 & 0 & 0 & 1 & 0 & 0 & 0 & 0  \\
    \end{block}
      \end{blockarray}\quad.
\]
Therefore we have
$$\dm_\cx\krn d^0_{\hn}=1,\qqq \dm_\cx\im d^0_{\hn}=3,$$
$$\dm_\cx\krn d^1_{H_{\mb{\scriptsize even}}}=5,\qqq \dm_\cx\krn d^1_{H_{\mb{\scriptsize odd}}}=4,$$
$$\dm_\cx\im d^1_{H_{\mb{\scriptsize even}}}=3,\qqq \dm_\cx\im d^1_{H_{\mb{\scriptsize odd}}}=4,$$
from which we conclude that
$$H^0(\W\sbu(H_{\mb{\scriptsize even}}))=\Z/2,\qq H^1(\W\sbu(H_{\mb{\scriptsize even}}))=(\Z/2)^2,\qq H^2(\W\sbu(H_{\mb{\scriptsize even}}))=\Z/2,$$
$$H^0(\W\sbu(H_{\mb{\scriptsize odd}}))=\Z/2,\qq H^1(\W\sbu(H_{\mb{\scriptsize odd}}))=\Z/2,\qq H^2(\W\sbu(H_{\mb{\scriptsize odd}}))=0.$$
Note that 
\beq{\krn d^0_{H_n}=H^0(\W\sbu(H_n))=\Z/2\ffall n\in\Z.\label{omg}}
Fern\'{a}ndez-Carmena \cite[3.4]{fer} showed that the Witt group of a complex surface is a birational invariant, so that
$$W(H_n)=W(\PP^2_\cx)=\Z/2\ffall n\in\Z.$$
On the other hand, since the rank is a local invariant, the localization map 
$$W(H_n)\ra W(H_{n, \eta})=\Gamma(\hn, \W^0(H_n)),$$
where $\eta\in\hn$ is the generic point, is an injection. By (\ref{omg}), we have
$$W(H_n)=H^0(\W\sbu(H_n)).$$
By the Purity Theorem \cite{OjP}, $\W\sbu$ is  a resolution of the sheaf $U\mt W(U)$ on $\hn$.

\appendix

\chapter{Technical lemmas}

We first note a useful lemma, which is easy to prove:

\lemm{Let $M,N,V$ be \amods, and
$$\phi:M\times N\ra V$$
an $A$-bilinear pairing. If
$$\ad\phi:M\ra \Hom_A(N,V),\qqq\ad\dg\phi:N\ra \Hom_A(M,V)$$
are the adjoints, then
\beq{\ad\phi=\Hom(\ad\dg\phi,V)\cc\rho_M,\qqq\ad\dg\phi=\Hom(\ad\phi,V)\cc\rho_N,
\label{magic}}
where
\bea{\Hom(\ad\phi,V):\Hom_A(M,V)&\la&\Hom_A(\Hom_A(N,V),V),\\
\Hom(\ad\dg\phi,V):\Hom_A(N,V)&\la&\Hom_A(\Hom_A(M,V),V),
} and 
\bea{\rho_M:M&\ra& \Hom_A(\Hom_A(M,V),V),\\
\rho_N:N&\ra& \Hom_A(\Hom_A(N,V),V)
}
are the canonical maps.
\label{first}
}

Note that if $\rho_M$ and $\rho_N$ are \isos, then $\ad\phi$ is bijective \ifof $\ad\dg\phi$ is bijective.

Now let $A:=\cx[x,y]$, $V^1_y:=\frac{\cx[x,y,1/y]}{\cx[x,y]}$, and denote $(-)\du:=\Hom_A(-,V^1_y)$. Pardon  proved the following :
\ben{\item If $M\in\cm^1_y(A)$, then $M\du\in\cm^1_y(A)$ \cite[1.13]{bigP}.\label{ap1}
\item $M\in\cm^1_y(A)$, then the canonical map
$$\rho_M:M\ra M\ddu$$\label{ap2}
is bijective \cite[1.17]{bigP}.
\item If
$$0\ra M'\ra M\ra M''\ra0$$
is an exact sequence in $\cm^1_y(A)$, then the induced sequence
$$0\la (M')\du\la M\du\la (M'')\du\la0$$
is exact \cite[1.6c]{P3}.\label{ap3}
\item If $M\in\cmya$ and $N\sst M$ is a submodule, then $N\in\cmya$ \cite[1.19]{bigP}.\label{ap4}
}

\lemm{Let $M,N\in\cmya$, and
$$\phi:M\times N\ra \vy$$
a bilinear pairing. Then the following are equivalent:
\ben{\item $\phi$ is \ns.
\item $\ad\phi$ is bijective.
\item $\ad\dg\phi$ is bijective.
\item $\ad\phi$ and $\ad\dg\phi$ are injective.
}
}

\pf{(1)\lra(2)\lra(3) follows from bijectivity of $\rho_M:M\ra M\ddu$. (1)\tra(4) is obvious.

(1)\tla(4) : Suppose that $\ad\phi$ and $\ad\dg\phi$ are injective.
Applying $\homa(-,\vy)$ to the short \esq
$$0\ra M\rf{\ad\phi}N\du\ra N\du/M\ra0$$
gives an injection
$$\cmya\owns N=N\ddu\hookleftarrow(N\du/M)\du.$$
Hence, $(N\du/M)\du\in\cmya$. Then $N\du/M=(N\du/M)\ddu\in\cmya$, so there is an \esq
$$0\la M\du\lf{(\ad\phi)\du}N\ddu\la(N\du/M)\du\la0,$$
so $(\ad\phi)\du$ is surjective.
Since $\ad\dg\phi=(\ad\phi)\du\cc\rho_N$ and $\rho_N$ is bijective, $\ad\dg\phi$ is surjective. A similar argument shows that $\ad\phi$ is surjecive.
}

\lemm{Let $M\in\cm^1_y(A)$, $N\sst M$ a submodule \st $M/N\in\cm^1_y(A)$, and
$$\phi:M\times M\ra \vy$$
a \ns\ symmetric $A$-bilinear form. Then $M/N\pp\in\cm^1_y(A)$, and $N=N\ppp$. Moreover, the induced pairings
$$\al:N\times M/N\pp\ra\vy,\qqq\be:N\pp\times M/N\ra V^1_y$$
are \ns.
\label{second}
}

\pf{Let
$$\ad\be:N\pp\ra\homa(M/N,V^1_y),\qqq\ad\dg\be:M/N\ra\homa(N\pp,V^1_y)$$
be the adjoints of $\be$. \Tcd
\xym{M \ar[r]_(.32)\sim^(.32){\ad\phi} & \homa(M,V^1_y) \\
N\pp \ar@{^(->}[u] \ar@{-->}[r]^(.3){\ad\be} & \ar@{^(->}[u] \homa(M/N,V^1_y)
}
Bijectivity of $\ad\phi$ implies that $\ad\be$ is bijective, so
$$\Hom(\ad\be,\vy):(N\pp)\du\la(M/N)\ddu$$
is bijective. By (\ref{magic}), $\ad\dg\be$ is then bijective. In particular, injectivity of $\ad\dg\be$ implies that $N\ppp\sst N$. Since $N\sst N\ppp$, we have $N=N\ppp$. 

Now
$$\ad\dg\al:M/N\pp\ra\homa(N,\vy)$$
is clearly injective, so $M/N\pp\in\cmya$. Hence, we can apply the same argument as above with $N$ replaced by $N\pp$ to conclude that $\ad\al$ is bijective.
}

The modules in $\cmya$ do not necessarily have finite length. However, if $M\in\cmya$ is $\cx[x]$-torsion-free, then one can define a notion similar to length; $y^{k-1}M/y^kM$ is a free module of finite rank over $\cx[x]$, so there is a finite chain of submodules
$$M=M^0\supset M^1\supset M^2\supset\cdots\supset M^n=0,$$
where $M^i/M^{i-1}\simeq\cx[x]$. We then define $\ell_y(M):=n$. Note the following:
\begin{itemize}
 \item $\elly$ is an additive function on $\cx[x]$-torsion-free modules in $\cmya$.
\item $M^i/M^{i-1}\in\cmya$, and it is $\cx[x]$-torsion-free.
\item If $M\in\cmya$ is $\cx[x]$-torsion-free, then so is $M\du$.
\end{itemize}

\lemm{Let $M\in\cmya$ be $\cx[x]$-torsion-free. If $N\sst M$ is a submodule, then $M/N\in\cmya$.
\label{slp}
}

\pf{It is clear that $M/N$ is $\cx[x]$-torsion-free. Hence,  $\dep_A(M/N)\geq1$, so $M/N\not\in\cmya$ implies that $M/N$ has dimension 2, i.e., it is of finite length. Then it is killed by a product of maximal ideals of the form $(x-a,y-b)\sst\cx[x,y]$, contradicting the $\cx[x]$-torsion-freeness.
}

\lemm{If $M\in\cmya$ is $\cx[x]$-torsion-free, then $\elly(M)=\elly(M\du)$.
\label{aft1}}

\pf{Let $\elly(M)=n$, so that there is a chain of submodules
$$M=M^0\supset M^1\supset M^2\supset\cdots\supset M^n=0$$
\st $M^i/M^{i-1}\simeq\cx[x]$. We prove by induction on $n$.
There is a short \esq\ of modules in $\cmya$ :
$$0\la M/M^1\la M\la M^1\la0.$$
Taking $\homa(-,\vy)$ gives an \esq\ \cite[1.6c]{P3}
$$0\ra(M/M^1)\du\ra M\du\ra(M^1)\du\ra0.$$
Note that as a $\cx[x,y]$-module, $M/M^1=\cx[x]$ is killed by $y$. Hence, the image of any \hmm\ $M/M^1\ra\vy$ lies in $(0:y)_{\vy}=\frac{\frac{1}{y}\cdot\cx[x,y]}{\cx[x,y]}\simeq\cx[x]$. Hence, $(M/M^1)\du\simeq\cx[x]$, so $\elly((M/M^1)\du)=1$. Then by the additivity of $\elly$ and the induction hypothesis,
$$\elly(M\du)=1+\elly((M^1)\du)=1+(n-1)=n=\elly(M).$$
}

\lemm{Let $M\in\cmya$ be $\cx[x]$-torsion-free, and
$$\phi:M\times M\ra\vy$$
a \ns\ symmetric bilinear form. If $N\sst M$ is a subspace, then 
$$\elly(N)+\elly(N\pp)=\elly(M).$$
\label{aft2}
}

\pf{By Lemma \ref{second}, the induced pairing
$$N\times M/N\pp\ra\vy$$
is \ns, so $N\simeq (M/N\pp)\du$. Hence, by Lemma \ref{aft1},
$$\elly(N)=\elly((M/N\pp)\du)=\elly(M/N\pp)=\elly(M)-\elly(N\pp).$$
}

\lemm{Let $M\in\cmya$ be $\cx[x]$-torsion-free,
$$\phi:M\times M\ra \vy$$
a \ns\ symmetric bilinear form, and $N\sst M$ is a totally isotropic subspace.
\ben{\item $2\cdot\elly(N)\leq\elly(M)$, and equality holds \ifof  $N$ is orthogonal, i.e., $N=N\pp$.
\item The  induced bilinear form
$$\bphi:\frac{N\pp}{N}\times\frac{N\pp}{N}\ra \vy$$
is \ns. Moreover, if $M$ has an orthogonal submodule, so does $N\pp/N$.
}
\label{final}
}

\pf{(1) : Since $N\sst N\pp$, $\elly(N)\leq\elly(N\pp)$. Hence, by Lemma \ref{aft2},
$$2\cdot\elly(N)\leq\elly(N)+\elly(N\pp)=\elly(M).$$
So $2\cdot\elly(N)=\elly(M)$ \ifof $\elly(N)=\elly(N\pp)$, i.e., if $N=N\pp$.

(2) : Since $N\pp/N\in\cmya$ by Lemma \ref{slp}, to prove nonsingularity of $\bphi$, it suffices to prove injectivity of $\ad\bphi$. But this is clear from $N\ppp=N$. 

Now suppose that $K\sst M$ is an orthogonal submodule. Let $\Lambda$ be the set of totally isotropic submodules containing $N$, and $S\in\Lambda$ a maximal element. We will show that $S\sst M$ is an orthogonal submodule. This implies that $S/N\sst N\pp/N$ is an orthogonal submodule, completing the proof.

First assume that $S\cap K=0$.
Since $S$ is totally isotropic, $2\cdot\elly(S)\leq\elly(M)$ by the first part of the proof. We will show that this is in fact an equality. Suppose, by way of contradiction, that $2\cdot \elly(S)<\elly(M)$. Since $K$ is orthogonal, $\elly(MWitt class)=2\cdot\elly(K)$ by the first part of the proof. Hence, $\elly(S)<\elly(K)$, and by Lemma \ref{aft2},
$$\elly(S\pp)=\elly(M)-\elly(S)>\elly(M)-\elly(K).$$
Hence, $\elly(M)<\elly(S\pp)+\elly(K)$, and this implies that $S\pp\cap K\neq0$, so there exists a non-zero element $m\in S\pp\cap K$. Since $S\cap K=0$ by assumption, $m\not\in S$. Since $m\in K$ and $K=K\pp$, $\psi(m,m)=0$. Hence, $S+(m)\sst M$ is a totally isotropic submodule strictly containing $S$, contradicting the maximality of $S$. Hence, $2\cdot\elly(S)=\elly(M)$, so $S\sst M$ is an orthogonal submodule by the first part of the proof.

Now for the general case, let $J:= S\cap K$. Then $J\sst M
$ is a totally isotropic submodule, so there is an induced bilinear form
$$\bpsi:J\pp/J\times J\pp/J\ra \vy.$$
We have shown that this is \ns. Note that
$$J\sst S\sst S\pp\sst J\pp,\qqq J\sst K= K\pp\sst J\pp.$$
$S/J\sst J\pp/J$ is maximal among totally isotropic submodules of $J\pp/J$. Moreover, 
$$(S/J)\cap(K/J)=0\sst J\pp/J,$$
so by the previous case, $S/J\sst J\pp/J$ is an orthogonal submodule. Hence, by Lemma \ref{second}, the induced pairings
$$\al:J\times M/J\pp\ra \vy,\qqq\be:S/J\times\frac{J\pp/J}{S/J}\ra \vy$$
are \ns. We will show that the pairing
$$\gam:S\times M/S\ra \vy$$
is \ns, which implies that $S=S\pp$.

There is an \esq
$$0\ra\frac{J\pp/J}{S/J}\ra\frac{M}{S}\ra\frac{M}{J\pp}\ra0.$$
Taking $\homa(-,\vy)$ gives a \cd
\xym{0  & \dssim_{\ad\be} S/J \ar[l] & \ar[d]_{\ad\gam} S \ar[l] & \dsim^{\ad\al} J \ar[l] & \ar[l] 0\\
0 & \ar[l] ((J\pp/J)/(S/J))\du & (M/S)\du \ar[l] & (M/J\pp)\du \ar[l] & \ar[l] 0
}
where the rows are exact. Hence, $\ad\gam$ is bijective, which implies that 
$$\ad\dg\gam:M/S\ra\homa(S,\vy)$$
is bijective. In particular, injectivity of $\ad\dg\gam$ implies that $S=S\pp$.
}

\rmk{From the proof of Proposition \ref{learn}, we know that there is an \iso\
\beq{1/y:W(\cx[x])\ras W(\cmya).\label{reall}}
Hence, every element of $W(\cmya)$ can be represented by a $\cx[x]$-torsion-free module. Lemma \ref{final}(2) and Lemma \ref{oba2} then suggests a way to obtain an inverse map of (\ref{reall}); let $[N,\psi]\in W(\cmya)$, where $N$ is $\cx[x]$-torsion-free, and $y^kN=0$ for some $k\geq1$. If $k=1$, then $N$ is a $\cx[x]$-module, and the image of $\psi$ lies in the image of the embedding
$$\V^0_{\cx[x]}=\cx[x]\os{1/y}{\hra}\frac{\cx[x,y,1/y]}{\cx[x,y]}=\V^1_{\cx[x,y],y}.$$
Hence, $[N,y\cdot\psi]\in W(\cx[x])$. If $k\geq2$, then $2k-2\geq k$, so $y^{k-1}N\sst N$ is a totally isotropic submodule, and there is an induced form
$$\bpsi:\frac{(y^{k-1}N)\pp}{y^{k-1}N}\times\frac{(y^{k-1}N)\pp}{y^{k-1}N}\ra\V^1_{\cx[x,y],y}.$$
By Lemma \ref{final}(2) and Lemma \ref{oba2}, $[N,\psi]=[\frac{(y^{k-1}N)\pp}{y^{k-1}N},\bpsi]$. Note that now we have $y^{k-1}\cdot\frac{(y^{k-1}N)\pp}{y^{k-1}N}=0$. Hence, by repeating this procedure, we will end up with the $k=1$ case above.
} 

\bibliographystyle{plain} 
\cleardoublepage
\normalbaselines 

\end{document}